\documentclass[a4paper,fleqn]{cas-sc}
\usepackage[square,numbers,sort,compress]{natbib}

\def\tsc#1{\csdef{#1}{\textsc{\lowercase{#1}}\xspace}}
\tsc{WGM}
\tsc{QE}
\usepackage{textcomp}
\usepackage{blindtext}
\hypersetup{
    colorlinks=true,
    linkcolor=blue,
    filecolor=magenta,      
    urlcolor=cyan,
}
\usepackage{color,soul}
\usepackage{placeins}

\usepackage{xcolor}
\usepackage{soul}
\usepackage[utf8]{inputenc}
\usepackage{array}
\usepackage[intoc, spanish]{nomencl}
\usepackage{subcaption}
\usepackage{textcomp}
\usepackage{graphicx}
\usepackage{amsthm}
\usepackage{amsmath,amssymb}
\usepackage{esvect}
\DeclareMathAlphabet{\mathbit}{OML}{cmm}{b}{it}
\usepackage{mathtools}
\usepackage{textcomp}
\usepackage{accents}

\makenomenclature
\usepackage{stackengine}
\usepackage{comment}
\usepackage{graphicx}
\usepackage[version=4]{mhchem}
\usepackage{siunitx}
\usepackage{longtable,tabularx}
\usepackage{mathtools}
\usepackage{pdfpages}
\usepackage{xcolor}
\usepackage{empheq}
\usepackage[most]{tcolorbox}
\usepackage{lscape}
\usepackage{amsfonts}

\tcbset{highlight math/.append style={boxrule=0pt,
                                      frame hidden,
                                      colback=yellow!40!white,
                                      sharp corners}}

\newtcbox{\mymath}[1][]{%
    nobeforeafter, math upper, tcbox raise base,
    enhanced, colframe=blue!30!black,
    colback=blue!30, boxrule=1pt,
    #1}
\usepackage[margins]{trackchanges}
\addeditor{JLM}
\addeditor{JH}
\begin{document} 

% \includepdf[pages=-]{RSER_Cover_Letter.pdf}
% \includepdf[pages=-]{RSER_Author_Checklist_Table.pdf}
%\includepdf[pages=-]{declaration_of_competing_interests.pdf}

\let\WriteBookmarks\relax
\def\floatpagepagefraction{1}
\def\textpagefraction{.001}

\shorttitle{Dynamic Modeling Of VSB WEC}    

\shortauthors{Shabara \& Abdelkhalik}  

%\title [mode = title]{Dynamic Modeling and Vibrations of Variable-Shape Wave Energy Converters} 
\title [mode = title]{Dynamic modeling of the motions of variable-shape wave energy converters} 
\tnotemark[1]
\tnotetext[1]{This document is the results of a research project funded by The National Science Foundation.}

\author[1]{Mohamed A. Shabara}[orcid = 0000-0002-8942-1412]
\ead{mshabara@iastate.edu}
\cormark[1]

\author[1]{Ossama Abdelkhalik}[orcid = 0000-0003-4850-6353]
\ead{ossama@iastate.edu}
\ead[url]{https://www.aere.iastate.edu/ossama/}

\affiliation[1]{organization={Department of Aerospace Engineering, Iowa State University},
            city={Ames},
            state={IA},
            postcode={50011}, 
            country={USA}}
            
\cortext[cor1]{Corresponding Author}

\begin{abstract}
In the recently introduced Variable-Shape heaving wave energy converters, the buoy changes its shape actively in response to changing incident waves. In this study, a Lagrangian approach for the dynamic modeling of a spherical Variable-Shape Wave Energy Converter is described. The classical bending theory is used to write the stress-strain equations for the flexible body using Love’s approximation. The elastic spherical shell is assumed to have an axisymmetric vibrational behavior. The Rayleigh-Ritz discretization method is adopted to find an approximate solution for the vibration model of the spherical shell. A novel equation of motion is presented that serves as a substitute for Cummins equation for flexible buoys. Also, novel hydrodynamic coefficients that account for the buoy mode shapes are proposed. The developed dynamic model is coupled with the open-source boundary element method software NEMOH. Two-way and one-way Fluid-Structure Interaction simulations are performed using MATLAB to study the effect of using a flexible shape buoy in the wave energy converter on its trajectory and power production. Finally, the variable shape buoy was able to harvest more energy for all the tested wave conditions.
%In the recently introduced Variable-Shape heaving wave energy converters, the buoy changes its shape actively in response to changing incident waves. In this study, a Lagrangian approach for the dynamic modeling for a spherical Variable-Shape Wave Energy Converter is described. The classical bending theory is used to write the stress-strain equations for the flexible body using Love’s approximation. The elastic spherical shell is assumed to have an axisymmetric vibrational behavior. The Rayleigh-Ritz discretization method is adopted to find an approximate solution for the vibration model of the spherical shell. A novel equation of motion is presented that serves as a substitute to Cummins equation for flexible buoys. Also, novel hydrodynamic coefficients that account for the buoy mode shapes are proposed. The developed dynamic model is coupled with the open source boundary element method software NEMOH. Two-way and one-way Fluid-Structure Interaction simulations are performed using MATLAB to study the effect of using a flexible shape buoy in the wave energy converter on its trajectory and power production. Finally, the variable shape buoy was able to harvest more energy for all the tested wave conditions.
\end{abstract}


\begin{comment}
\begin{highlights}
\item Dynamic Modeling of Variable-shape Buoy Wave Energy Converters
\item A Lagrangian approach for writing the equations of motion using a Rayleigh-Ritz approximation
\item Comparison of harvested power from Variable-shape Buoy and Fixed-Shape Buoy Wave Energy Converters
\item Using damping control force, the Variable-shape Buoy Wave Energy Converter behaves like a Fixed-shape Buoy Wave Energy Converter with a PD control
\item Demonstrating the advantage of using  Variable-shape Buoy Wave Energy Converters to reduce the need for reactive power
%\item Rayleigh-Ritz approximation is used to obtain approximate solution for the displacement vector.
%\item Flexible buoys would leverage the waves and behave like a rigid buoy that has reactive power; indeed the reactive power in this case is obtained from the waves themselves.
%\item Kinetic and potential energies of asymmetric free vibrating spherical shells. 
\end{highlights}
\end{comment}

% Keywords


\maketitle

\section*{Nomenclature}
{\renewcommand\arraystretch{1.0}
\noindent\begin{longtable*}{@{}l @{\quad=\quad} l@{}}
$c$  & Damping Coefficient (Ns/m)\\
$E$  & Young Modulus of Elasticity (MPa)\\
$h$  & Sphere thickness (m)\\
$V$  & Volume (m$^3$)\\
$t$  & Time (sec)\\
$r$  & Sphere Radius (m)\\
$\rho$ & Material Destiny (kg/m$^3$)\\
$\nu$  & Poisson's Ratio\\
$\eta$  & Rayleigh-Ritz Coefficient\\
$\omega$  & Frequency (rad/sec)\\
$\boldsymbol{1}$ & Identity matrix $=[\boldsymbol{1}_1 \; \boldsymbol{1}_2 \; \boldsymbol{1}_3]$\\ 
\multicolumn{2}{@{}l}{Subscripts and Superscripts}\\
hydro & hydrodynamic\\
$b$ & buoyancy \\
hst & hydrostatic \\
$pto$  & Power Take-off Unit \\
$w$  & Water \\
ref & reference \\
rad & radiation \\
$\dagger$ & Pseudo Inverse 
\end{longtable*}}

\section{Introduction}

Oceans are colossal reservoirs of energy of particularly high density energy \cite{clement2002wave}. The total theoretical ocean energy potential is estimated to be $29.5$ PWh/yr \cite{energyreport2016}, which is more than the US electric power needs in $2020$. Despite its significant potential, ocean energy is still a very small portion of the overall renewable energy production \cite{RENEWABLES2017}.

One widely used concept for harvesting wave power is the heaving Wave Energy Converter (WEC). In its simplest form, this point absorber device may consist of a floating buoy connected to a vertical hydraulic cylinder (spar) attached at the bottom of the seabed.  
As the buoy moves due to the wave and control forces, the hydraulic cylinders drive hydraulic motors, and the motors drive a generator \cite{Li2012392}.
The forces on the floating buoy are the excitation, radiation, and hydrostatic forces \cite{Bacelli2014}. The excitation force is due to the wave field and the buoy's geometry. The motion of the buoy itself creates waves which in turn create the radiation forces. The hydrostatic force accounts for the buoyancy force and weight of the buoy. 
In most of the current wave harvesting devices, the WEC has a Fixed-Shape Buoy (FSB).
The equation of motion for a heave-only 1-DoF FSB WEC is \cite{cummins1962impulse}:
\begin{equation}\label{eq:Cummins}
m\ddot{z}(t) = \overbrace{\int_{-\infty}^{\infty}h_f(\tau)\eta(t - \tau,z)d\tau}^\text{excitation force $f_e$} + f_s  \overbrace{-\mu\ddot{z}(t) - \int_{-\infty}^{t}h_r(\tau)\dot{z}(t - \tau)d\tau}^\text{radiation force $f_r$} - u
\end{equation}where $m$ is the buoy mass, $u$ is the control force, $t$ is the time, $z$ is the heave displacement of the buoy from the sea surface, and $f_s$ is the hydrostatic force that reflects the spring-like effect of the fluid. The $\eta$ is the wave surface elevation at buoy centroid, $f_e$ is the excitation force, and $h_f$ is the impulse response function defining the excitation force in heave. The radiation force is $f_r$, where $\mu$ is a frequency-dependent added mass, and $h_r$ is the impulse response function defining the radiation force in heave.

For a FSB, the convolution integral part of the radiation force $f_r$ in Eq.~\eqref{eq:Cummins} can be approximated using a state space model of $N$ states, $\vec{x}_r = [x_{r1}, \cdots , x_{rN}]^T$, which outputs the radiation force \cite{Yu1995}:
\begin{eqnarray}\label{eqXr}
\dot{\vec{x}}_r  =  A_r \vec{x}_r + B_r \dot{z} & \text{and} & f_r  =  C_r \vec{x}_r,
\end{eqnarray}where the constant radiation matrices $A_r$ and $B_r$ are obtained by approximating the impulse response function in the Laplace domain, as detailed in several references such as \cite{Perez2009a}.

Two important aspects impact the energy converted from oceans: the control force and  the buoy's shape. For an FSB WEC, linear dynamic models are widely used in control design e.g. \cite{Bacelli2011,Falnes2007,6803979,Ringwood2014IEEE,Scruggs20131,Li2012392,Hals2011_2,abdulkadir4274235optimal,9867173,9705961}.
Usually, the control is designed to maximize the mechanical power of the WEC; many references adopt different approaches to achieve optimality. Often, the resulting control forces have a spring-like component in addition to the resisting force \cite{Fusco2011ReacCtrl,WECSingArc}. Hence, the Power Take-Off (PTO) unit needs to have a bidirectional power flow capability, which is typically complex and more expensive. 
In analyzing the shapes of FSB WECs, the use of any non-cylindrical shape requires the use of non-linear hydro models \cite{hamada2021numerical,hamada2022numerical}. This is the reason that most studies assume cylindrical shape of the FSB WEC. 

From an economic perspective, the cost of having the complex bidirectional power flow PTO to maximize harvested energy is high. Moreover, the structure of a FSB WEC needs to be designed to withstand very high loads at peak times despite operating at a much less load most of the time. This impacts the structural design and increases the cost. To mitigate this peak load, geometry controlled OSWEC was recently proposed in references \cite{tom2015preliminary,tom2016spectral,tom2016development,kelly2017development}, where controllable surfaces, along with a wave-to-wave control, are used to maximize power capture, increase capacity factor, and reduce design loads. The latter controlled-geometry OSWEC changes shape only when the wave climate changes, and hence it can be considered similar to the case of an FSB WEC when it is not in the transition from one geometry to another. %More discussion of this type is presented in \MySection~\ref{IMSHMV}.
A geometry control of the overtopping WEC is proposed in reference \cite{victor2011effects}. The slope angle and crest freeboard of the device is made adaptive to the sea conditions by geometry control.  Reference \cite{kurniawan2012characteristics} proposed a variable flap angle pitching device.  The resonance characteristics of the WEC can be altered by controlling the angle of the flap.  Later, reference \cite{kurniawan2017wave} proposed a floating airbag WEC that has a longer resonance period without implementing phase control.

%The wave energy converters can be classified either based on their position (shoreline, near-shore, and offshore devices) or principle of operation or based on the power take-off (PTO) technique. Drew et al. \cite{drew2009review} classified the WEC based on the principle of operation into three main categories as follows

%\begin{enumerate}
%\item Attenuators, such as the Pelamis which are devices that lie parallel to the predominant wave direction.
%\item Point absorbers, which have small dimensions relative to the incident wave wavelength.  
%\item Terminators, such as the Salter's Duck which are devices that lie perpendicular to the predominant wave directions.
%\end{enumerate}

%\hl{WEC can be further classified based on their shell flexibility; fixed shape "FSB WEC" and flexible membrane WECs "FlexWEC"}. A typical optimal control requires reactive power, which needs to be supplied by a complex and expensive Power Takeoff (PTO) unit. the energy produced by an optimal control is also significantly lower than expected (derived without considering PTO dynamics) and the resulting poor power quality (significant fluctuation) \cite{zou2020modeling,zou2021numerical}. A passive control does not need reactive power though the performance is not optimal.

The concept of a variable-shape buoy (VSB) WEC was recently introduced to reduce the complexity of the PTO. A VSB WEC changes its shape continuously. The wave/WEC interaction produced by the VSB WEC can be leveraged to produce more power without adding complexity to the PTO unit.
%The difference between a VSB WEC and a variable geometry WEC is the rate of change of the surface of the WEC shape relative to the incident wave forces, i.e., the variable geometry WEC is an FSB WEC that changes its geometry occasionally, on the other hand, the VSB WEC changes their shapes continuously. The wave/WEC interaction produced by the VSB WEC can be tuned to produce more power without complexity accompanied by the reactive power PTO units.
Specifically, Zou et al. \cite{zou2020modeling} proposed the Variable-Shape point absorber; their original design comprises a pressurized gas chamber attached to a set of multiple controllable moving panels. This VSB WEC is controlled by a simple linear damping PTO unit \cite{zou2020modeling}. A low-fidelity dynamic model is derived in \cite{zou2020modeling} to demonstrate the superiority of the VSB WEC compared to the FSB WEC. The average power harvested using the VSB WEC in \cite{zou2020modeling} is about 18\% more compared to the FSB WEC. %Also, it is noticed that the velocity of the VSB WEC are higher than the FSB WEC. 

In another study, references \cite{zou2020numerical,ZouJMSE2021NumSimVSBWEC} present three-dimensional two-way Fluid-Structure Interaction (FSI) high fidelity simulations, using the ANSYS software package, to simulate a spherical VSB WEC. The WEC in \cite{zou2020numerical,ZouJMSE2021NumSimVSBWEC} has a hyper-elastic hollow shell of radius $2$ m. The internal volume contains trapped gas that helps in creating a restoring moment. The device is simulated in a numerical wave tank (NWT) with dimensions $80 \times 60 \times 60$ m$^3$, and damping regions at the sides and at the outlet of the NWT. The free surface height was at $40$ m. Their simulation captured the highly non-linear behavior of the VSB WEC and showed an enhancement in the heave displacement and velocity for the VSB WEC compared to a similar-size FSB WEC.
Reference \cite{shabara2021numerical} presents a study that uses a similar approach but applies a passive control force. The PTO force is dependent on the WEC heave velocity ($F_{PTO} = -c \dot{x}$), where $c$ is a constant damping coefficient. The results showed an increase in the heave displacement and velocity for the VSB WEC over the similar-size FSB WEC. The results also show an increase in the harvested energy of about $8\%$. 

As can be seen from the above discussion, studies on VSB WECs currently use either high fidelity numerical software tools or rough approximate low fidelity tools \cite{ogden2021review} for simulations. The high fidelity tools are computationally expensive, and the above low fidelity simulations cannot capture important features in this FSI phenomenon.% since rigid moving panels are used to mimic the flexible body dynamics. %Also, although WEC-Sim software is capable of modeling flexible bodies with generalized body modes, however the flexibility should should be added as a degree of freedom in the fluid solver not in the structure domain which limits the .  

In this work, Lagrangian mechanics and Rayleigh-Ritz approximation are used to derive novel equation of motions for variable shape wave energy converters $\left(\text{Eqs.} \eqref{eq:105_3}, \eqref{eq:105_4}\, \text{and} \, \eqref{EoMF_VSB}\right)$. Novel expressions for the generalized added mass, damping, hydrostatic, and excitation forces and coefficients are also derived analytically for two-way and one-way FSI schemes (sections \ref{sec:const_forced} and \ref{sec:one-two_way}). The paper is organized as follows: In section \ref{sec:math}, the kinematics of the flexible buoy are presented, and then the kinetic and potential energies for the special case of spherical buoys are presented \cite{shabara2022Dynamics}, noting that the use of spherical buoys is only to demonstrate the utility of the model, and the derived aforementioned generalized hydrodynamic/hydrostatic coefficients and forces can be used with any flexible buoy geometry. In section \ref{sec:EOM}, the equations of motion for flexible shell buoys are derived using Lagrangian mechanics for the free unconstrained vibration. In section \ref{sec:const_forced}, a generalized form for the hydrodynamic forces and coefficients are proposed, for regular and irregular waves. Section \ref{sec:FSI} presents the most general form of the proposed equation of motion for VSB WECs and discusses the two-way FSI scheme for the proposed model. To reduce the computational time associated with the two-way FSI schemes, Reynolds averaging is applied to obtain novel hydrodynamic coefficients for the one-way FSI in section \ref{sec:one-two_way}. The model validation is discussed in section \ref{sec:valid}. Finally, the numerical simulation results for the one-way and two-way FSI for regular and irregular waves are presented in section \ref{sec:results}.

%In this work the body frame rotational movements are omitted and heave only motion is assumed, the dynamic model is coupled with the boundary element method (BEM) open source code Nemoh; A comprehensive comparison between one-way and two-way FSI simulations was performed. Also, Nemoh was modified to run in parallel when simulating irregular waves which significantly reduces the computational time by orders of magnitude. 

%\section{Reference Frames and Notations}
\section{Kinetic and Potential Energies of Spherical Shell Buoys} \label{sec:math}

The derivation of the equation of motion in the current study uses Lagrangian mechanics; hence, the calculation of the kinetic and potential energies is required. This section starts with the layout of the used reference frames, then a description of the domain discretization technique used in the current study. The system's kinematics are derived in subsection \ref{subsec:Kinematics}; the result is then used to calculate the kinetic and potential energies in subsections \ref{subsec:KE} and \ref{subsec:PE}, respectively.

\begin{figure}
\centering
\includegraphics[scale=0.45]{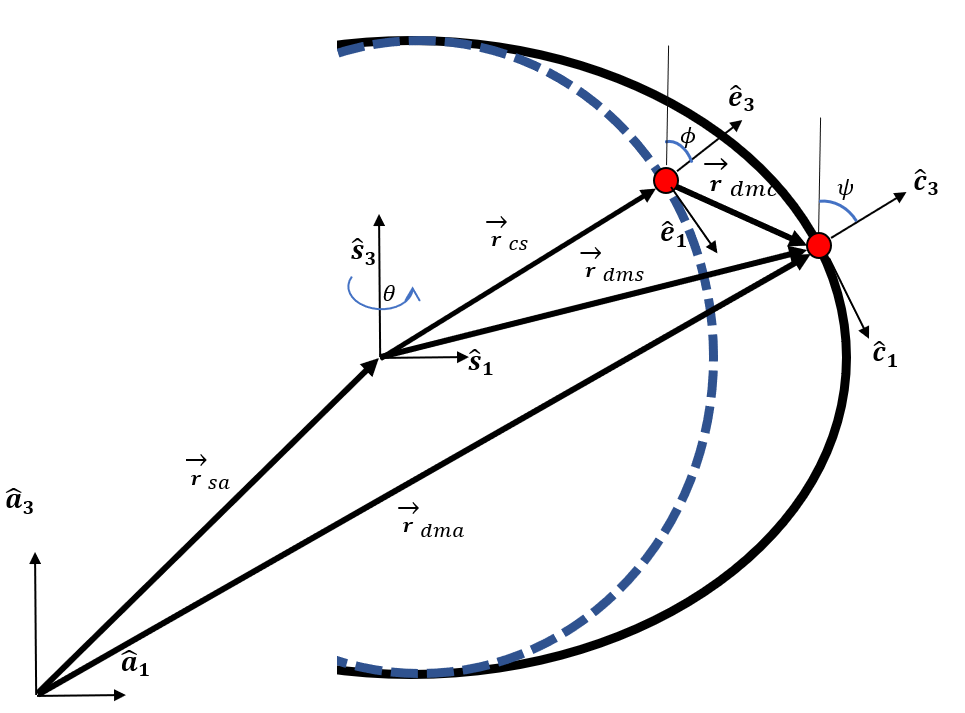}
\caption{Deformed (Solid Black Line) and non-deformed Buoy (Dashed Blue Line)}
\label{fig:1}
\end{figure}

Consider a flexible buoy for which the non-deformed shape is spherical. As shown in Fig. (\ref{fig:1}), the inertial frame is denoted as $\hat{\boldsymbol{a}}$ and can be described as:
\begin{equation}
\hat{\boldsymbol{a}} = \left[ \hat{\boldsymbol{a}}_1 , \hat{\boldsymbol{a}}_2 , \hat{\boldsymbol{a}}_3 \right ] 
%=\mbox{diag}\{\begin{matrix} %\hat{a}_1 & \hat{a}_2 & %\hat{a}_3 \end{matrix} \} 
\end{equation}
%where diag is the diagonal function.
The body-fixed frame $\hat{\boldsymbol{s}}$ is attached to the buoy's center of mass. 
Any point on the buoy's surface can be specified using the two coordinates $\phi$ and $\theta$, as illustrated in Fig. (\ref{fig:1}).
Consider an infinitesimal mass at the surface of the buoy; we introduce the reference frame $\hat{\boldsymbol{e}}$ which is attached to that infinitesimal mass on the buoy's surface before deformation, and its third axis $\hat{\boldsymbol{e}}_3$ is aligned with the radius of the non-deformed buoy shape.

Hence, the reference frame $\hat{\boldsymbol{e}}$ is obtained by rotating $\hat{\boldsymbol{s}}$ by an angle $\theta$ around the $\hat{\boldsymbol{s}}_3$ then by an angle $\phi$ around the second intermediate frame as follows. 
\begin{equation}
    C_{es}(\phi,\theta) = C_2(\phi)C_3(\theta) 
\end{equation}where $C_i(x)$ represents a fundamental transformation matrix of a single rotation of angle $x$ about the coordinate $i$, $i=2,3$.
The reference frame $\hat{\boldsymbol{c}}$ is centered at the infinitesimal mass on the surface such that $\hat{\boldsymbol{c}}_3$ is normal to the surface.
The angle $\psi$ is the angle between $\hat{\boldsymbol{s}}_3$ and $\hat{\boldsymbol{c}}_3$.
In the analysis presented in this paper, it is assumed that the deformations are axisymmetric about the $\hat{\boldsymbol{a}}_3$ axis; hence the axis $\hat{\boldsymbol{c}}_2$ is always perpendicular to the page. The axes $\hat{\boldsymbol{a}}_2$, $\hat{\boldsymbol{s}}_2$, and $\hat{\boldsymbol{e}}_2$ are also perpendicular to the page.
If the shape is not deformed from its original spherical shape, then the frames $\hat{\boldsymbol{e}}$ and $\hat{\boldsymbol{c}}$ coincide. The frames $\hat{\boldsymbol{e}}$ and $\hat{\boldsymbol{c}}$ become different, in general, when the shape is deformed. 
For the FSB WEC the reference frames $\hat{\boldsymbol{e}}$ and $\hat{\boldsymbol{c}}$ coincide; this applies to the VSB WEC at the initial time before deformation.

References \cite{zou2021numerical,shabara2021numerical,zou2020numerical} carried out high-fidelity FSI simulations using ANSYS for spherical VSB WECs, and they found that the steady-state response for the VSB WECs is close to being axisymmetric; thus in this work the VSB WEC response is assumed to be axisymmetric to simplify the analysis.

The coordinate transformation matrix from the $\hat{\boldsymbol{a}}$ frame to the $\hat{\boldsymbol{s}}$ frame is computed in this paper using the 3-2-1 Euler angle sequence as $C_{sa}(\alpha,\beta,\gamma) =  C_1(\alpha)C_2(\beta) C_3(\gamma)$. %\hl{omitted part}

Since the changes in the buoy shape are assumed axisymmetric, we can express the deformation vector (displacement) as a function of only the angle $\phi$ and the time $t$. This deformation vector can be expressed in the $\hat{\boldsymbol{e}}$ frame as: \begin{align} \label{eq:1}
 \vec{\boldsymbol{r}}_{dmc}(\phi,t) = \begin{bmatrix} u (\phi,t) & 0 & v (\phi,t) \end{bmatrix}{}^T \end{align} where the second component (normal to the page plane) is set to zero because of the axisymmetry of the deformation, $u(\phi,t)$ is the displacement component in the $\hat{\boldsymbol{e}}_1$ direction, and $v(\phi,t)$ is the displacement component in the $\hat{\boldsymbol{e}}_3$ direction. In this paper, each of these displacement components is assumed a series of separable functions; that is, each term in their series can be expressed as a product of two functions, one of them depends only on $\phi$ and the other depends only on $t$. Moreover, the Rayleigh-Ritz approximation is used to obtain an approximate solution for the displacement vector as discussed in the following section. 

\subsection{Rayleigh-Ritz Approximation}
The approximation method used in this work is the Rayleigh-Ritz method.
Each component of the displacement vector $\vec{\boldsymbol{r}}_{dmc}$ is assumed to have the following form \cite{hogan2015dynamic,shabara2022Dynamics,9867800}:
\begin{equation}
u(\phi,t) = \sum_{n=1}^N \Psi_{n}^{\phi}(\phi)\eta_n(t) = \underbrace{[\Psi_{1}^{\phi} \ldots \Psi_{N}^{\phi }]}_{\mathbf{\Psi}_{e}^{\phi}}  \underbrace{\begin{bmatrix}
    \eta_{1}(t) \\
    \vdots \\
    \eta_{N}(t) 
\end{bmatrix}}_{\boldsymbol{\eta}(t)}= \mathbf{\Psi}_{e}^{\phi}(\phi) \boldsymbol{\eta}(t) \label{eq:51_5}
\end{equation}

\begin{equation}
v(\phi,t) = \sum_{n=1}^N \Psi_{n}^{r}(\phi)\eta_n(t)= \underbrace{[\Psi_{1}^{r} \ldots \Psi_{N}^{r}]}_{\mathbf{\Psi}_{e}^{r}} 
\underbrace{\begin{bmatrix}
    \eta_{1}(t) \\
    \vdots \\
    \eta_{N}(t) 
\end{bmatrix}}_{\boldsymbol{\eta}(t)}= \mathbf{\Psi}_{e}^{r}(\phi) \boldsymbol{\eta}(t) \label{eq:52_5}
\end{equation}where the functions $\Psi_{n}^{\phi}$ and $\Psi_{n}^{r}$ are trial (admissible) functions of $\phi$ and the functions $\eta_n$ are functions of time $t$, $\forall$ $n=1, \cdots, N$.
Therefore, the displacement vector can be expressed in the $\hat{\boldsymbol{e}}$ frame as follows:
\begin{equation}\label{eq:5}
\vec{\boldsymbol{r}}_{dmc}(\phi,t) = \underbrace{\begin{bmatrix}
\mathbf{\Psi}_{e}^{\phi}(\phi) \\
\mathbf{0} \\
\mathbf{\Psi}_{e}^{r}(\phi) 
\end{bmatrix}}_{\boldsymbol{\Phi}_e} \boldsymbol{\eta}(t) = \boldsymbol{\Phi}_e(\phi) \boldsymbol{\eta}(t) \end{equation}

For a spherical shape buoy, the Legendre functions of the first kind $P_n$ \cite{kreyszig2009advanced,shabara2022Dynamics} can serve as shape functions for the Ritz-Rayleigh method to satisfy the essential geometrical (Dirichlet) boundary conditions  \cite{naghdi1962vibrations,hogan2015dynamic,soedel2005vibrations,newman1994wave,raouf1990non} as follows:

\begin{equation} \label{eq:6}
\Psi_{n}^{\phi}(\phi) = A \frac{d P_n (\cos(\phi))}{d \phi} \text{, and} \:\:  \Psi_{n}^{r}(\phi) = A \frac{(1+(1+\nu))\Omega_n^2}{1-\Omega_n^2}   P_n (\cos(\phi))
\end{equation} where the coefficients of the equations above form an eigenvector for the Legendre differential equation, i.e. the constant "A" can take any real value. $\Omega_n^2$ is a dimensionless frequency parameter expressed as \cite{soedel2005vibrations,shabara2022Dynamics,9867800}:
\begin{equation} \label{eq:7}
\Omega^2_{n} = \frac{1}{2(1-\nu^2)} (\bar{A} \pm \sqrt{\bar{A}^2 -4m\bar{B}})
\end{equation}where $\nu$ is the Poisson's ratio, and
\begin{align}
m &= n(n+1)-2,  \: n \in \mathbb{Z}^{+} \label{eq:10}\\
\bar{B} &= 1+ \nu^2 +\frac{1}{12} [(m+1)^2 - \nu^2] \label{eq:9}\\
\bar{A} &= 3(1+\nu) + m+\frac{1}{2} \left[ \frac{h}{r} \right]^2 (m+3)(m+1+\nu) \label{eq:8}
\end{align}
From \cite{soedel2005vibrations,naghdi1962vibrations,9867800} the natural frequencies in radians per second for spherical shells are calculated using Eq. (\ref{eq:11})
\begin{equation} \label{eq:11}
\omega_n^2 = \frac{E}{r^2 \rho \Omega_n^2}
\end{equation} where $E$ is Young's Modulus, $r$ is the non-deformed radius of the shell, and $\rho$ is the density of the shell material.
When $n = 0$, the vibration mode corresponds to the breathing mode (volumetric or pulsating modes) which is a pure radial vibration mode \cite{soedel2005vibrations,nayfeh2006axisymmetric,raouf1990non}. For $n > 0$, the $\pm$ sign in Eq.~\eqref{eq:7} yields the modes corresponding to the membrane vibration modes and bending vibration modes. The bending vibration modes are obtained when using the negative sign; these modes are sensitive to the $h/r$ ratio. On the other hand, the membrane modes are insensitive to the change in the $h/r$ ratio. Due to the extensional motion of the buoy, only the membrane vibration modes are used in the current work. To obtain the approximated equations of motion using the Rayleigh-Ritz method, the approximated displacement vector needs to be substituted in the kinetic and strain energy equations as follows.

\subsection{Kinematics of a Flexible Spherical Buoy -  Free Vibration} \label{subsec:Kinematics}
This subsection is concerned with calculating the ${}^a\dot{\vec{\boldsymbol{r}}}_{dma}$ vector as it is crucial for the calculation of the kinetic energy of the spherical shell due to the translation and rotational motions as well as the deformation of the sphere external shell.
As shown in Fig. (\ref{fig:1}), the position vector of a point on the surface of the deformed sphere in the inertial frame "$\hat{\boldsymbol{a}}$" is expressed as \cite{hogan2015dynamic,shabara2022Dynamics}: \begin{equation} \label{eq:24}
\vec{\boldsymbol{r}}_{dma}=\vec{\boldsymbol{r}}_{sa} +\vec{\boldsymbol{r}}_{cs} +\vec{\boldsymbol{r}}_{dmc}
\end{equation}
The velocity vector is expressed as
\begin{equation}
{}^a\dot{\vec{\boldsymbol{r}}}_{dma} ={}^a\dot{\vec{\boldsymbol{r}}}_{sa} +{}^a\dot{\vec{\boldsymbol{r}}}_{cs} +{}^a\dot{\vec{\boldsymbol{r}}}_{dmc} 
\end{equation}
Note that the left superscript denotes the reference frame used to describe the vector. Accordingly, ${}^a\dot{\vec{\boldsymbol{r}}}_{dma}$ for a buoy constraint from rotations can be expressed as \cite{shabara2022Dynamics}:
\begin{align}
{}^a\dot{\vec{\boldsymbol{r}}}_{dma} & = \begin{bmatrix}
C_{sa} & C_{se} 
\end{bmatrix} \underbrace{\begin{bmatrix}
{}^a\dot{\vec{\boldsymbol{r}}}_{sa} \\
{}^e\dot{\vec{\boldsymbol{r}}}_{dmc}
\end{bmatrix}}_{\dot{\vec{\mathbf{x}}}} \label{eq:31_3}
\end{align}

\subsection{Kinetic Energy for a Flexible Spherical Buoy } \label{subsec:KE} 
The total kinetic energy of the buoy is expressed in Eq. (\ref{eq:32_3}).
\begin{equation} \label{eq:32_3}
\mathcal{T} = \frac{1}{2} \int_S   {}^a\dot{\vec{\boldsymbol{r}}}_{dma} \cdot {}^a\dot{\vec{\boldsymbol{r}}}_{dma} dm \\
\end{equation}where $S$ denotes the surface of the buoy, Substituting Eq.~\eqref{eq:31_3} in Eq. (\ref{eq:32_3}) to get \cite{shabara2022Dynamics},
\begin{align}
 \mathcal{T}  =  \underbrace{ \frac{1}{2}  {}^a\dot{\vec{\boldsymbol{r}}}_{sa}^T \boldsymbol{m} {}^a\dot{\vec{\boldsymbol{r}}}_{sa}}_{\mathcal{T}_{x}}  +
\underbrace{ \frac{1}{2}\int_s  {}^e\dot{\vec{\boldsymbol{r}}}_{dmc}^T   {}^e\dot{\vec{\boldsymbol{r}}}_{dmc}  dm}_{\mathcal{T}_{\boldsymbol{s}}} 
\label{eq:47_4}
\end{align}
where $\boldsymbol{m} = \boldsymbol{1}m$ is the VSB mass matrix. The terms ${\mathcal{T}_{x}}$ and ${\mathcal{T}_{\boldsymbol{s}}}$ are the kinetic energies associated with the transnational motion and shell deformations, respectively.

The kinetic energy $\mathcal{T}_{\boldsymbol{s}}$ is approximated by substituting Eq.s \eqref{eq:51_5} and \eqref{eq:52_5} into Eq. \eqref{eq:47_4}; to get \cite{hogan2015dynamic,shabara2022Dynamics}: 
\begin{align} 
\mathcal{T}_{\boldsymbol{s}} &= \frac{1}{2} \dot{\boldsymbol{\eta}}^T \underbrace{\Bigg\{ 2 \pi \rho h \int^{\pi}_{0} \left( \mathbf{\Psi}_{e}^{{\phi}^T} \mathbf{\Psi}_{e}^{\phi} +  \mathbf{\Psi}_{e}^{{r}^T} \mathbf{\Psi}_{e}^{r}\right)  r^2 \sin \phi d \phi \Bigg\}}_{\boldsymbol{M}_{ee}} \dot{\boldsymbol{\eta}} = \frac{1}{2} \dot{\boldsymbol{\eta}}^{T} \boldsymbol{M}_{ee} \dot{\boldsymbol{\eta}}  \label{eq:14}
\end{align}

\subsection{Potential Energy for Flexible Spherical Buoys} \label{subsec:PE}
The strain energy-displacement expressions (membrane strains) for axisymmetric shells can be written as \cite{forbes2010dynamic,rao2007vibration,soedel2005vibrations}:
\begin{align}
\varepsilon_{\phi\phi} & = \frac{1}{r} \left(  \frac{\partial u}{\partial\phi}+v  \right) \label{eq:15.1} \\ 
\varepsilon_{\theta\theta} & = \frac{1}{r} \left( u \cot(\phi) + v \right) \label{eq:16.1}
\end{align}
The total strain energy can be found in References \cite{shabara2022Dynamics,rao2007vibration}. Combining the sphere elastic strain and the gravitational potential energies in Eqs. (\ref{eq:19}) and (\ref{eq:29}) yield the total potential energy of the spherical shell buoy \cite{hogan2015dynamic,shabara2022Dynamics}:
\begin{align}
\mathcal{\pi} =  \mathcal{U}_s + \mathcal{G} &=   \frac{1}{2} \frac{Eh}{1-\nu^2} \int_0^{2 \pi} \int_0^\pi \left\lbrace{ \varepsilon_e^{\phi \phi}}^2+ {\varepsilon_e^{\theta \theta}}^2 + 2 \nu {\varepsilon_e^{\phi\phi}} \varepsilon_e^{\theta \theta} \right\rbrace  r^2 \sin(\phi) d \theta d \phi + mg  \boldsymbol{1}_3^T \vec{\boldsymbol{r}}_{sa} \label{eq:18}
\end{align}   
The strain energy is approximated by first substituting Eqs. \eqref{eq:51_5} and \eqref{eq:52_5} in Eqs. \eqref{eq:15.1} and \eqref{eq:16.1} to get:
\begin{align}
\varepsilon_{\phi\phi} & =  \frac{1}{r} \left(  \frac{\partial \Psi_e^\phi}{\partial\phi}+\Psi_e^r\right)\boldsymbol{\eta} \label{eq:15} \\ 
\varepsilon_{\theta\theta} & =   \frac{1}{r} \left( \Psi_e^\phi \cot(\phi) + \Psi^r_e \right) \boldsymbol{\eta} \label{eq:16}
\end{align}
Then the strain energy equation in Eq.~\eqref{eq:18} becomes \cite{hogan2015dynamic,shabara2022Dynamics,forbes2010dynamic}:
\begin{align}
\mathcal{U}_{\boldsymbol{s}} & = \frac{1}{2}{\boldsymbol{\eta}}^T \Bigg[ \frac{2\pi Eh}{1-\nu^2} \int_0^{\pi} \Bigg\{ \left(    \frac{\partial \Psi_e^\phi}{\partial\phi}+\Psi_e^r\right)^T \left(  \frac{\partial \Psi_e^\phi}{\partial\phi}+\Psi_e^r\right) +  \left( \Psi_e^\phi \cot(\phi) + \Psi^r_e \right)^T  \left( \Psi_e^\phi \cot(\phi) + \Psi^r_e \right) \nonumber \\
&+\nu \left(  \left(  \frac{\partial \Psi_e^\phi}{\partial\phi}+\Psi_e^r\right)^T   \left( \Psi_e^\phi \cot(\phi) + \Psi^r_e \right) +  \left( \Psi_e^\phi \cot(\phi) + \Psi^r_e \right)^T    \left(  \frac{\partial \Psi_e^\phi}{\partial\phi} + \Psi_e^r\right) \right)\Bigg\}\sin(\phi)  d\phi \Bigg] \boldsymbol{\eta} \\
&= \frac{1}{2} {\boldsymbol{\eta}}^T \boldsymbol{K}_{ee} \boldsymbol{\eta} \label{eq:19}
\end{align}
where 

\begin{align}
\boldsymbol{K}_{ee}   &= \frac{2\pi Eh}{1-\nu^2} \int_0^{\pi} \left\lbrace \zeta_{\phi\phi}^T \zeta_{\phi\phi} + {\zeta_{\theta \theta}}^T \zeta_{\theta \theta} +\nu \left( \zeta_{\phi\phi}^T \zeta_{\phi\phi} +{\zeta_{\theta \theta}}^T \zeta_{\theta \theta}\right) \right\rbrace  \sin(\phi) d\phi \\
\zeta_{\phi\phi} &=   \left(  \frac{\partial \Psi_e^\phi}{\partial\phi}+\Psi_e^r\right) \label{eq:20}\\ 
\zeta_{\theta \theta} &= \left( \Psi_e^\phi \cot(\phi) + \Psi^r_e \right) \label{eq:21}
\end{align}

%\subsubsection{Gravitational Energy}
Considering the buoy as a spherical shell, the gravitational energy can be expressed as \cite{hogan2015dynamic,shabara2022Dynamics}:
\begin{equation}
\mathcal{G}  = \int_s -{\vec{\boldsymbol{g}}} \cdot \vec{\boldsymbol{r}}_{sa} dm = mg  {\hat{\boldsymbol{a}}_3}^T \vec{\boldsymbol{r}}_{sa}   \label{eq:29}
\end{equation}
where $\vec{\boldsymbol{g}} = -g \hat{\boldsymbol{a}}_3$.

\section{Unconstrained Equations of Motion For Flexible shell Buoys} \label{sec:EOM}

The unconstrained equations of motion are here derived as a first step towards writing the constrained equations of motion. The Lagrangian for this buoy system can be written as the summation of three quantities \cite{junkins2009analytical}:
\begin{align}
    \mathcal{L} =\mathcal{L}_D + \mathcal{L}_B + \int_s \hat{\mathcal{L}} d \phi
\end{align}
where $\mathcal{L}_D(t,\vec{\mathbf{x}},\dot{\vec{\mathbf{x}}})$ is associated with the discrete coordinates, $\hat{\mathcal{L}}$ is the Lagrangian density function and it is a function of the discrete and distributed parameter coordinates. Using the Rayleigh-Ritz approximations transfers the distributed parameter terms to discrete parameter terms, i.e., in this case $\hat{\mathcal{L}}= \boldsymbol{0}$, $\mathcal{L}_B$ is associated with the boundaries.
In this work there are no boundary terms in the Lagrangian equation, i.e., $\mathcal{L}_B = \boldsymbol{0}$. Now, we will derive the equations of motion related to the discrete coordinates.

% \subsection{Equation of Motion Associated with Discrete Coordinates}
The Lagrangian for the discrete coordinates is expressed as:
\begin{align} 
\mathcal{L}_D = \mathcal{T}_D - \mathcal{\pi}_ D
            =  \frac{1}{2}   {}^a\dot{\vec{\boldsymbol{r}}}_{sa}^T \boldsymbol{m} {}^a\dot{\vec{\boldsymbol{r}}}_{sa} + \frac{1}{2} \dot{\boldsymbol{\eta}}^T \boldsymbol{M}_{ee} \dot{\boldsymbol{\eta}} - \frac{1}{2} {\boldsymbol{\eta}}^T \boldsymbol{K}_{ee} \boldsymbol{\eta}  -mg  \boldsymbol{1}_3^T \vec{\boldsymbol{r}}_{sa} \label{eq:30_5}
\end{align}

The Lagrange Equation for the discrete coordinates is expressed as:
\begin{equation} \label{eq:49_3}
\frac{d}{dt}\left(\frac{\partial \mathcal{L}_D}{\partial \dot{\vec{\mathbf{x}}}}\right) - \frac{\partial \mathcal{L}_D}{\partial{\vec{\mathbf{x}}}}= \mathbf{0}
\end{equation}

To write the equations of motion of the discrete coordinates, we first write:
\begin{align}
&\frac{\partial \mathcal{L}_D}{\partial{\vec{\boldsymbol{r}}_{sa}}} = - mg\boldsymbol{1}_3 \text{, and}\;\;\; \frac{d}{dt} \left( \frac{\partial \mathcal{L}_D }{\partial \dot{\vec{\boldsymbol{r}}}_{sa}} \right)   =  \boldsymbol{m} {}^a\ddot{\vec{\boldsymbol{r}}}_{sa} \\
&   \frac{\partial \mathcal{L}_s}{\partial{\boldsymbol{\eta}}}  =- \boldsymbol{K}_{ee} \boldsymbol{\eta} \\
& \frac{\partial \mathcal{L}_s}{\partial\dot{\boldsymbol{\eta}}}  = \boldsymbol{M}_{ee} \dot{\boldsymbol{\eta}}   \text{, and} \, \, \frac{d}{dt}\left(\frac{\partial \mathcal{L}_s}{\partial\dot{\boldsymbol{\eta}}}\right) =  \boldsymbol{M}_{ee} \ddot{\boldsymbol{\eta}}
\end{align} 
Noting that, unlike \cite{shabara2022Dynamics}, the generalized mass matrix is constant due to the no C.G rotations assumption. The equations of motion for the translational motion and the buoy deformations are:
\begin{align}
     \boldsymbol{m} {}^a\ddot{\vec{\boldsymbol{r}}}_{sa} + mg\boldsymbol{1}_3 &=\boldsymbol{0}\label{eq:34_5} \\
     \boldsymbol{M}_{ee} \ddot{\boldsymbol{\eta}} + \boldsymbol{K}_{ee} \boldsymbol{\eta} &= \boldsymbol{0}\label{eq:26_1} 
\end{align} The equations of motion described by Eq. \eqref{eq:34_5} and \eqref{eq:26_1} can be extended to include damping coefficient matrices as follows: \begin{align}
     \boldsymbol{m} {}^a\ddot{\vec{\boldsymbol{r}}}_{sa} + \boldsymbol{D}_x {}^a\dot{\vec{\boldsymbol{r}}}_{sa} + mg\boldsymbol{1}_3 &= \boldsymbol{0} \label{eq:36_5}\\
\boldsymbol{M}_{ee} \ddot{\boldsymbol{\eta}} + \boldsymbol{D}_{ee} \dot{\boldsymbol{\eta}} + \boldsymbol{K}_{ee} \boldsymbol{\eta} &= \boldsymbol{0}  \label{eq:36_6}
\end{align}

where $\boldsymbol{D}_x$ is the damping metric for the translation motions and $\boldsymbol{D}_{ee}$ is a {proportional} damping matrix that is assumed to be a function of the mass and stiffness matrices as follows:
\begin{equation} \label{eq:23_1}
\boldsymbol{D}_{ee} = \alpha_d \boldsymbol{M}_{ee} + \beta_d \boldsymbol{K}_{ee} \nonumber
\end{equation}  
where the $\alpha_d$ and $\beta_d$ are real scalars called the mass and stiffness matrix multipliers with units 1/sec and sec, respectively \cite{alipour2008study,liu1995formulation,thompson2017ansys}. Noting that the modal damping, in that case, is dependant only on the rate of change of flexibility states, the addition of damping terms here is important to damp out the transient vibration response and obtain the steady state response within a finite settling time.

Combining the equations of motion from Eqs. \eqref{eq:36_5} and \eqref{eq:36_6} yields the unforced equation of motion of a flexible buoy constraint from rotations \cite{hogan2015dynamic,shabara2022Dynamics}:
\begin{align}
& \boldsymbol M {\ddot{\vec{\mathbf{x}}}}  + \boldsymbol{D} {\dot{\vec{\mathbf{x}}}}+ \begin{bmatrix}
mg \boldsymbol{1}_3^T & (\boldsymbol{K}_{ee} \boldsymbol{\eta})^T
\end{bmatrix}^T = \boldsymbol{0} \label{eq:50_2}
\end{align} 
where
\begin{align}
\dot{\vec{\mathbf{x}}} &= \begin{bmatrix}
{}^a\dot{\vec{\boldsymbol{r}}}_{sa}^T   &
\dot{\boldsymbol{\eta}}^T  \end{bmatrix}^T_{(3+N) \times 1}    \\
\boldsymbol M  &=  \text{diag} \{ \begin{matrix} 
       \boldsymbol{m}& \boldsymbol{M}_{ee} \end{matrix} \}_{(3+N) \times (3+N)} \\
\boldsymbol{D} &= \text{diag} \{ \begin{matrix} 
     \boldsymbol{D}_x & \boldsymbol{D}_{ee} \end{matrix} \}_{(3+N) \times (3+N)}  
\end{align}

\section{Forced Constrained Equations of Motion} \label{sec:const_forced}
In this section, a Lagrangian mechanics approach is used to derive a new generalized form for  the hydrodynamic forces. Then, novel expressions for the generalized hydrodynamic coefficients are derived. The Lagrange equations for the forced motion of discrete coordinates are expressed as \cite{junkins2009analytical}:
\begin{align} 
\frac{d}{dt}\left(\frac{\partial \mathcal{L}}{\partial \dot{\mathbf{x}}}\right) - \frac{\partial \mathcal{L}}{\partial{\mathbf{x}}} &= \boldsymbol{Q} \nonumber  \\ & =  \boldsymbol{Q}^{c} +\boldsymbol{Q}^{pto} + \boldsymbol{Q}^{ext} + \boldsymbol{Q}^{rad} +\boldsymbol{Q}^{\infty}+ \boldsymbol{Q}^{hst} \label{eq:49_5}
\end{align}
where, $\boldsymbol Q$ is a column vector of generalized forces and can be expressed as the summation of the generalized constraint force ``$\boldsymbol Q^{c}$" \cite{shabara2022Dynamics}, generalized PTO force ``$\boldsymbol Q^{pto}$", generalized excitation force ``$\boldsymbol Q^{ext}$", generalized interia ``$\boldsymbol Q^{\infty}$", generalized radiation force $\boldsymbol Q^{rad}$, and generalized hydrostatic force ``$\boldsymbol Q^{hst}$". Noting that the application of the constraint forces on the VSB WEC was demonstrated in Ref. \cite{shabara2022Dynamics}.

To compute the generalized force, the general transformation takes the form \cite{junkins2009analytical}: 
\begin{equation} \label{eq:60_1g} 
  \mathbit{Q}_j =  \sum_{i=1}^{3+N} f_i \cdot \frac{\partial \vec{ r_i}}{\partial q_j}
\end{equation} 

The next subsections assumes that the buoy is constrained from surging and swaying and only heave motion is allowed.

\subsection{Power Take-off Unit Force} 
The PTO can be either active or passive \cite{osti_10292551}. In this paper, the PTO is assumed to be passive damping. To apply the PTO effect for that case, one can either include its effect in the $\boldsymbol{D}_x$ matrix in Eq.~\eqref{eq:50_2} or as an external force in Eq. \eqref{eq:73_3} \cite{shabara2022Dynamics}. \begin{equation} \label{eq:73_3}
    \vec{\mathbit{f}}^{pto} =   -c \: {}^a\dot{\vec{\mathit{\boldsymbol{r}}}}_{sa,3} \: \boldsymbol{a}_3
\end{equation} where $\vec{\mathbit{f}}^{pto}$ is the damping force and $c$ is the damping coefficient.

Hence, the generalized force corresponding to the PTO force takes the form: \begin{equation} \label{eq:60_1} 
  \mathbit{Q}_j^{pto} =  \vec{\boldsymbol{f}}^{pto} \cdot \frac{\partial \vec{r}_{sa}}{\partial q_j}, \, \, \forall j=1, 2, \cdots 1+N
\end{equation}

We can then write the PTO generalized forces using Eq. \eqref{eq:60_1} and Eq. \eqref{eq:73_3} as follows:
\begin{align}
 \mathbit{Q}_1^{pto} = \vec{\boldsymbol{f}}^{pto}   \cdot \frac{\partial  \vec{\boldsymbol{r}}_{sa} }{\partial \vec{\boldsymbol{r}}_{sa,3}} &= -c {}^a\dot{\vec{\mathit{\boldsymbol{r}}}}_{sa,3}  \boldsymbol{a}_3 \cdot \boldsymbol{1}_3 = -c {}^a\dot{\vec{\mathit{\boldsymbol{r}}}}_{sa,3}  \\
 \mathbit{Q}_2^{pto} = \vec{\boldsymbol{f}}^{pto} \cdot \frac{\partial  \vec{\boldsymbol{r}}_{sa} }{\partial {\eta_1} } &= 0, \:\:\:
\hdots \hdots \:\:\:
 \mathbit{Q}_{1+N}^{pto} =  \vec{\boldsymbol{f}}^{pto} \cdot \frac{\partial  \vec{\boldsymbol{r}}_{sa} }{\partial {\eta_N} } = 0
\end{align} 
In a compact form, 
\begin{equation}
    \mathbit{Q}^{pto} = - \begin{bmatrix}
      c \: {}^a\dot{\vec{\mathit{\boldsymbol{r}}}}_{sa,3} \:  \\ \boldsymbol{0} 
    \end{bmatrix}
\end{equation}  Note that the equation above accounts for the heave-only motion (the rows related to surging and swaying were omitted). 

\subsection{Generalized hydrodynamic Forces For Regular Waves} \label{subsec:3.1}
In this section, the hydrodynamic forces for the VSB WEC are derived in the generalized coordinates, including the generalized added mass, the generalized damping coefficients, the generalized hydrostatic coefficient, and the generalized excitation force. The frequency domain excitation, radiation, and hydrostatic pressures are assumed to be known or calculated using BEM solvers (ex, NEMOH, WAMIT, AQWA, Capytaine). These BEM solvers discretize the buoy shell into several mesh elements/panels, and the hydro pressures on each can be extracted from the solver.  

% then extend the derivation to derive the generalized forces for irregular waves. 
 
\subsubsection{Derivation of the Generalized Excitation  Force $\boldsymbol{Q}^{ext}$:} \label{sssec:exh}

The generalized hydrodynamic excitation force "$\mathbit{Q}^{ext}$" is computed by considering the submerged volume of a buoy, and the hydro excitation force acting on the $i^{th}$ mesh element is expressed as 
\begin{align}
\vec{f}_i^{ext} &= - {F_{ext}}_i   {\hat{\boldsymbol{c}}_3}_i \label{eq:27_4} \\
F_{ext,i} &= p_{ext,i}A_i  
\end{align} where $A_i$ and ${\hat{\boldsymbol{c}}}_{3,i}$ are the $i^{th}$ element surface area and mesh element normal, respectively, and $p_{ext,i}$ is the time domian excitation pressure generated by the waves on the $i^{th}$ mesh element. Note that the $t$ subscript denotes that these variables are in the time domain.

The time domain excitation pressure $p_{ext,i}$ is calculated using Eq.~\eqref{eq:29_1} \cite{zou2020modeling} \begin{equation} \label{eq:29_1}
 p_{ext,i} =  \Re \left( p_{ex,i} \eta e^{\hat{i} \omega t} \right)
\end{equation}     

where $\eta$ is the wave amplitude ($\eta = H/2$) \cite{wacher2010mathematical}, $ H$ is the wave height, $\omega$ is the wave frequency, $p_{ex,i}$ is the frequency domain excitation pressure coefficient, and $\hat{i}$ equals $\sqrt{-1}$. 

The frequency domain excitation pressure coefficient is the summation of the Froude-Krylov (KF) $p_{\text{kf},i}$ and the diffraction $p_{\text{df},i}$ pressure coefficients. The Froude-Krylov (KF) pressure is generated by the unsteady pressure field due to the undisturbed incident wave, while the diffraction pressure is caused due to the WEC disturbing the waves \cite{zou2020modeling}. \begin{equation}
p_{ex,i} = p_{\text{kf},i} +p_{\text{df},i} \: \: \: \in \mathbb{C}
\end{equation}

The $p_{\text{df},i}$ is obtained numerically using the BEM open source tools, while the value for Froude-Krylov (KF) is calculated analytically depending on the sea depth \cite{falnes2002ocean}. For deep sea conditions the KF pressure can be expressed in the following notation: \begin{equation} \label{eq:87}
p_{\text{kf},i} = \rho_w g e^{k[({\vec{\boldsymbol{r}}_{dms3}}(\phi_w)- {\vec{\boldsymbol{r}}_{dms}}_3 )+i{{\vec{\boldsymbol{r}}_{dms}}_1}_i]}
\end{equation}
 
where $k = \omega^2/g$ is the wave number \cite{zou2020modeling}, and ${\vec{\boldsymbol{r}}_{dms3}}(\phi_w)$ is the vertical component of the position vector of the mesh element coinciding with the free surface. 

From Eq.~\eqref{eq:60_1g}, the generalized hydro force can be written in the following form:
\begin{equation} \label{eq:60_1z} 
  \mathbit{Q}_j^{ext} = \vec{\boldsymbol{f}}^{ext} \cdot \frac{\partial \vec{r}_{dma}}{\partial q_j} , \, j = 1,2, ..., 1+N
\end{equation} 

Substituting Eq.~\eqref{eq:27_4} into Eq.~\eqref{eq:60_1z} to get:\begin{align}
\mathbit{Q}_{1,i}^{ext} = \vec{\boldsymbol{f}}^{ext}_i   \cdot \frac{\partial  \vec{\boldsymbol{r}}_{dma} }{\partial \vec{\boldsymbol{r}}_{sa3}} &= - {F_{exh}}_i {\hat{\boldsymbol{c}}_3}_i \cdot \boldsymbol{1}_3 =  - {F_{ext}}_i \boldsymbol{1}_3^T {\hat{\boldsymbol{c}}_3}_i =- {F_{ext}}_i \cos{\psi_i}    \\
 \mathbit{Q}_{2,i}^{ext} = \vec{\boldsymbol{f}}^{ext}_i \cdot \frac{\partial  \vec{\boldsymbol{r}}_{dma} }{\partial {\eta_1} } &= - {F_{ext}}_i {\hat{\boldsymbol{c}}_3}_i  \cdot \mathbf{\Phi}_{e,i}(:,1) , \:\:\:
\hdots, \:\:\:
 \mathbit{Q}_{1+N,i}^{ext} =  \vec{\boldsymbol{f}}^{ext}_i \cdot \frac{\partial  \vec{\boldsymbol{r}}_{sa} }{\partial {\eta_N} } = - {F_{ext}}_i {\hat{\boldsymbol{c}}_3}_i  \cdot \mathbf{\Phi}_{e,i}(:,N) 
\end{align} 

Therefore the generalized excitation force for heave only motion is as follows:, \begin{align}
    \mathbit{Q}^{ext} &=- \sum_{i=1}^{m} F_{ext,i} \begin{bmatrix}  \cos{\psi_i} \\
     {\hat{\boldsymbol{c}}_3}_i  \cdot \mathbf{\Phi}_{e,i}(:,1) \\ 
    \vdots  \\
     {\hat{\boldsymbol{c}}_3}_i  \cdot \mathbf{\Phi}_{e,i}(:,N)  \\
    \end{bmatrix} =- \sum_{i=1}^{m} p_{ext,i} A_i  \begin{bmatrix}  \cos{\psi_i} \\
    {\hat{\boldsymbol{c}}_3}_i  \cdot \mathbf{\Phi}_{e,i}(:,1) \\ 
    \vdots  \\
    {\hat{\boldsymbol{c}}_3}_i  \cdot \mathbf{\Phi}_{e,i}(:,N)  \\
    \end{bmatrix}   \label{eq:81}
\end{align} where $m$ is the total number of mesh panels. Accordingly, the proposed generalized excitation force coefficient for VSB WEC is expressed as: \begin{equation}
    \boldsymbol{Ex}(\omega,t) =- \sum_{i=1}^m p_{ex,i}(\omega,t) A_i  \begin{bmatrix}  \cos{\psi_i} \\
    {\hat{\boldsymbol{c}}_3}_i  \cdot \mathbf{\Phi}_{e,i}(:,1) \\ 
    \vdots  \\
    {\hat{\boldsymbol{c}}_3}_i  \cdot \mathbf{\Phi}_{e,i}(:,N)  \\
    \end{bmatrix} \: \: \: \in \mathbb{C} \label{eq:53_3}
\end{equation} where $\boldsymbol{Ex}$ is a column vector. Note that for the two-way FSI simulations, each time step requires calculating a new $\boldsymbol{Ex}$ due to its dependence of the buoy's geometry. Finally, Eq.~\eqref{eq:81} can be reformulated as: \begin{equation}
    \boldsymbol{Q}^{ext}(t) = \Re \left( \boldsymbol{Ex}(\omega,t) \eta e^{\hat{i} \omega t} \right)
\end{equation}

\subsubsection{Derivation of the generalized hydrostatic Force $\boldsymbol{Q}^{hst}$:} \label{sec:buoyant}
The derivation of the generalized hydrostatic force $\boldsymbol{Q}^{hst}$ will yield an expression for the generalized hydrostatic matrix as demonstrated in this subsection. The hydrostatic pressure on the $i^{th}$ mesh element is computed as \cite{falnes2002ocean}: \begin{equation}
    p_{hst,i} =-\rho_w g (h_{\text{ref}}- {\vec{\boldsymbol{r}}_{dma}}_{3,i} ) = -\rho_w g z_i  
\end{equation} $h_\text{ref}$ is the water free-water level, and ${\vec{\boldsymbol{r}}_{dma}}_3$ is expressed as: \begin{align}
\vec{\boldsymbol{r}}_{dma} &=\vec{\boldsymbol{r}}_{sa} +  r {\hat{\boldsymbol{e}}_3} +{\boldsymbol{\Phi}_e} \boldsymbol{\eta}\\
{\vec{\boldsymbol{r}}_{dma}}_{3,i} &= \begin{bmatrix}
 C_{sa}(3,3) & [C_{se}\boldsymbol{\Phi}_e(\phi_i)](3,:) 
\end{bmatrix} \vec{\mathbf{x}} \label{eq:65_3}
\end{align} The hydrostatic resorting force is the difference between the buoy weight and the hydrostatic forces and can be calculated as follows: \begin{align}
F_{hst,i} &= p_{hst,i}A_i  =\rho_w g  A_i  (h_\text{ref} -  {\vec{\boldsymbol{r}}_{dma3}}_i ) =\rho_w g  A_i  z_i  \label{eq:27_1} \\
\vec{f}^{hst}_{tot}  &= -  \sum_{i=1}^{N} \left( {F_{hst}}_i   {\hat{\boldsymbol{c}}_3}_i - \delta m_i \,  g {\hat{\boldsymbol{a}}_3} \right) = -  \sum_{i=1}^{N} \left( {F_{hst}}_i  {C_{ac}}_i - \delta m_i \,  g  \right){\hat{\boldsymbol{a}}_3}  \\
&= - \rho_w g \sum_{i=1}^{N} \left(   A_i  z_i {C_{ac}}_i{\hat{\boldsymbol{a}}_3}  - \delta m_i \, {\hat{\boldsymbol{a}}_3}   \right)\label{eq:27_3}
\end{align} where $\delta m_i$ is the mass of the $i^{th}$ mesh element, at equilibrium, the buoy is half submerged, and the weight forces balances out the hydrostatic forces; therefore, Eq.~\eqref{eq:27_1} can be expressed as: \begin{align}
\vec{\boldsymbol{f}}^{hst} &= - \rho_w g \sum_{i=1}^{N} \left(  A_i  \Delta z_i {\hat{\boldsymbol{c}}_{3,i}}     \right) =- \rho_w g \sum_{i=1}^{N} \left(   A_i  (z_{1,i} - z_i) {\hat{\boldsymbol{c}}_{3,i}}   \right)  \\
&=- \rho_w g \sum_{i=1}^{N} \left(   A_i  (h_\text{ref} -  {\vec{\boldsymbol{r}}_{dma3}}_{i,0}  - h_\text{ref} +  {\vec{\boldsymbol{r}}_{dma3}}_i ) {\hat{\boldsymbol{c}}_{3,i}}   \right)  \\
&= - \rho_w g \sum_{i=1}^{N} \left(   A_i( {\vec{\boldsymbol{r}}_{dma3}}_{i,0}   -  {\vec{\boldsymbol{r}}_{dma3}}_i ) {\hat{\boldsymbol{c}}_{3,i}}   \right) \label{eq:41_1}
\end{align} Note that at the initial time the VSB WEC is not deformed ${\vec{\boldsymbol{r}}_{dmc}} (t=0)={\vec{\boldsymbol{r}}_{dmc}}_0 = \boldsymbol{0}$. Assuming that the inertial frame is on the free-water surface ${\vec{\boldsymbol{r}}_{sa}}_0 = \boldsymbol{0}$, we can write:
\begin{align} 
{\vec{\boldsymbol{r}}_{dma3}}_{i,0}   -  {\vec{\boldsymbol{r}}_{dma3}}_i &= {\vec{\boldsymbol{r}}_{sa3}}_{0} +{\vec{\boldsymbol{r}}_{cs3}}_{i,0} +{\vec{\boldsymbol{r}}_{dmc3}}_{i,0} - {\vec{\boldsymbol{r}}_{sa3}}_{i} -{\vec{\boldsymbol{r}}_{cs3}}_{i} - {\vec{\boldsymbol{r}}_{dmc3}}_{i} \nonumber \\
&=  - \left( {\vec{\boldsymbol{r}}_{sa3}}  + {\vec{\boldsymbol{r}}_{dmc3}}_{i} \right)
\end{align} Accordingly, Eq~\eqref{eq:41_1} can be expressed as: \begin{equation} \label{eq:64_2}
    \vec{\boldsymbol{f}}^{hst} = \rho_w g \sum_{i=1}^{N} \left(   A_i  \left( {\vec{\boldsymbol{r}}_{sa3}}  + {\vec{\boldsymbol{r}}_{dmc3}}_{i} \right) {\hat{\boldsymbol{c}}_{3,i}}   \right)\\ 
\end{equation} From Eq.~\eqref{eq:60_1g}, the generalized hydrodynamic stiffness can be expressed as: \begin{equation} \label{eq:60_1y} 
  \mathbit{Q}_j^{hst} = \vec{\boldsymbol{f}}^{hst} \cdot \frac{\partial \vec{r}_{dma}}{\partial q_j} , \, j = 1,2, ..., 1+N
\end{equation} Substituting Eq.~\eqref{eq:64_2} into Eq.~\eqref{eq:60_1y} to get: \begin{align}
\mathbit{Q}_{1,i}^{hst} &= \vec{\boldsymbol{f}}^{hst}    \cdot \frac{\partial  \vec{\boldsymbol{r}}_{dma} }{\partial \vec{\boldsymbol{r}}_{sa3}} =  \rho_w g A_i   \left( {\vec{\boldsymbol{r}}_{sa3}}  + {\vec{\boldsymbol{r}}_{dmc3}}_{i} \right) {\hat{\boldsymbol{c}}_3}_i \cdot \boldsymbol{1}_3 =  \rho_w g    A_i    \left( {\vec{\boldsymbol{r}}_{sa3}}  + {\vec{\boldsymbol{r}}_{dmc3}}_{i} \right)  \cos{\psi_i}    \\
 \mathbit{Q}_{2,i}^{hst} &= \vec{\boldsymbol{f}}^{hst}_i \cdot \frac{\partial  \vec{\boldsymbol{r}}_{dma} }{\partial {\eta_1} } = \rho_w g A_i   \left( {\vec{\boldsymbol{r}}_{sa3}}  + {\vec{\boldsymbol{r}}_{dmc3}}_{i} \right)  {\hat{\boldsymbol{c}}_3}_i  \cdot \mathbf{\Phi}_{e,i}(:,1) \\
 \mathbit{Q}_{1+N,i}^{hst} &=  \vec{\boldsymbol{f}}^{hst}_i \cdot \frac{\partial  \vec{\boldsymbol{r}}_{sa} }{\partial {\eta_N} } = \rho_w g A_i   \left( {\vec{\boldsymbol{r}}_{sa3}}  + {\vec{\boldsymbol{r}}_{dmc3}}_{i} \right)  {\hat{\boldsymbol{c}}_3}_i  \cdot \mathbf{\Phi}_{e,i}(:,N) 
\end{align} 

Therefore, \begin{align}
    \mathbit{Q}^{hst} &= \sum_{i=1}^{m}  \rho_w g A_i    \begin{bmatrix}  \cos{\psi_i} \\
     {\hat{\boldsymbol{c}}_3}_i  \cdot \mathbf{\Phi}_{e,i}(:,1) \\ 
    \vdots  \\
     {\hat{\boldsymbol{c}}_3}_i  \cdot \mathbf{\Phi}_{e,i}(:,N)  \\
    \end{bmatrix} \left( {\vec{\boldsymbol{r}}_{sa3}}  + {\vec{\boldsymbol{r}}_{dmc3}}_{i} \right) \label{eq:78_3} \end{align} Subsitute Eq.~\eqref{eq:78_3} into Eq.~\eqref{eq:65_3} to get: \begin{align}
\mathbit{Q}^{hst} &=  \underbrace{\sum_{i=1}^{m}  \rho_w g A_i  \begin{bmatrix}  \cos{\psi_i} \\
     {\hat{\boldsymbol{c}}_3}_i  \cdot \mathbf{\Phi}_{e,i}(:,1) \\ 
    \vdots  \\
     {\hat{\boldsymbol{c}}_3}_i  \cdot \mathbf{\Phi}_{e,i}(:,N)  \\
    \end{bmatrix} \begin{bmatrix}
 C_{sa}(3,3) & [C_{se}\boldsymbol{\Phi}_e(\phi_i)] (3,:) 
\end{bmatrix}}_{-\boldsymbol{K}_h} \vec{\mathbf{x}}
\end{align}

Accordingly, the generalized hydrostatic stiffness coeffcient matrix for heave-only motion is expressed as:

\begin{equation} \label{eq:48}
 \boldsymbol{K}_h = - \sum_{i=1}^{m}  \rho_w g A_i  \begin{bmatrix}  \cos{\psi_i} \\
     {\hat{\boldsymbol{c}}_3}_i  \cdot \mathbf{\Phi}_{e,i}(:,1) \\ 
    \vdots  \\
     {\hat{\boldsymbol{c}}_3}_i  \cdot \mathbf{\Phi}_{e,i}(:,N)  \\
    \end{bmatrix} \begin{bmatrix}
 1 & [C_{se}\boldsymbol{\Phi}_e(\phi_i)] (3,:) 
\end{bmatrix} 
\end{equation} noting that $ \boldsymbol{K}_h$ is a full matrix which makes the equation of motion statically coupled. The derived matrix is also non-symmetric.

\subsubsection{Derivation of the Generalized Radiation Force $\boldsymbol{Q}^{rad}$:}
The derivation of the generalized radiation force $\boldsymbol{Q}^{rad}$ will yield an expression for the generalized radiation damping matrix, and this will be demonstrated in this subsection.
The portion of the radiation force that is dependent on the buoy velocity for regular waves is expressed as:  
\begin{align}
\vec{f}_i^{rad} &= - {F_{rad}}_i   {\hat{\boldsymbol{c}}_3}_i \label{eq:27} \\
F_{rad,i} &= p_{rdt,i}A_i   
\end{align} where $A_i$ and ${\hat{\boldsymbol{c}}}_{3,i}$ are the $i^{th}$ element surface area and it's normal, respectively. 

Time domain radiation pressure for only heave motion is calculated as follows \cite{zou2020modeling}: \begin{equation}
p_{rdt,i} = p_{rdc,i} {\dot{\vec{\boldsymbol{r}}}_{dma}}_3
\end{equation} where $p_{rdc,i}$ is the real part of the radiation damping pressure coefficient \cite{zou2020modeling}, and it is calculated using Eq.~\eqref{eq:31}; also, $\dot{\boldsymbol{r}}_{dma}$ for motion can be calculated as using Eq.~\eqref{eq:32}: \begin{align}
p_{rdc,i} &= \Re \left( p_{rd,i} \right) \label{eq:31} \\
{}^a\dot{\vec{\boldsymbol{r}}}_{dma} &= \begin{bmatrix}
C_{sa} & C_{se} \boldsymbol{\Phi}_e
\end{bmatrix} \dot{\vec{\mathbf{x}}} \label{eq:32}
\end{align} where $p_{rd,i} \in \mathbb{C} $ is the radiation pressure coefficient \cite{zou2020modeling} obtained from the BEM solver, and the superscript $a$ denotes that the vector is expressed in the $a$ frame. 

Substituting Eq.~\eqref{eq:27} into Eq.~\eqref{eq:60_1z} to get:\begin{align}
\mathbit{Q}_{1,i}^{rad} &= \vec{\boldsymbol{f}}^{rad}_i   \cdot \frac{\partial  \vec{\boldsymbol{r}}_{dma} }{\partial \vec{\boldsymbol{r}}_{sa,3}} = - {F_{rad}}_i {\hat{\boldsymbol{c}}_3}_i \cdot \boldsymbol{1}_3 =  - {F_{rad}}_i \boldsymbol{1}_3^T {\hat{\boldsymbol{c}}_3}_i =- p_{rdc,i} A_i \dot{\boldsymbol{r}}_{dma3}  \cos{\psi_i}   \\
 \mathbit{Q}_{2,i}^{rad} &= \vec{\boldsymbol{f}}^{rad}_i \cdot \frac{\partial  \vec{\boldsymbol{r}}_{dma} }{\partial {\eta_1} } = - {F_{rad}}_i {\hat{\boldsymbol{c}}_3}_i  \cdot \mathbf{\Phi}_{e,i}(:,1) = - p_{rdc,i} A_i \dot{\boldsymbol{r}}_{dma3} {\hat{\boldsymbol{c}}_3}_i  \cdot \mathbf{\Phi}_{e,i}(:,1) \\
 \mathbit{Q}_{1+N,i}^{rad} &=  \vec{\boldsymbol{f}}^{rad}_i \cdot \frac{\partial  \vec{\boldsymbol{r}}_{sa} }{\partial {\eta_N} } = - {F_{rad}}_i {\hat{\boldsymbol{c}}_3}_i  \cdot \mathbf{\Phi}_{e,i}(:,N)  = - p_{rdc,i} A_i \dot{\boldsymbol{r}}_{dma3} {\hat{\boldsymbol{c}}_3}_i  \cdot \mathbf{\Phi}_{e,i}(:,N) 
\end{align} 

Therefore, \begin{align}
    \mathbit{Q}^{rad} &=- \sum_{i=1}^{m} {F_{rad}}_i \begin{bmatrix}  \cos{\psi_i} \\
     {\hat{\boldsymbol{c}}_3}_i  \cdot \mathbf{\Phi}_{e,i}(:,1) \\ 
    \vdots  \\
     {\hat{\boldsymbol{c}}_3}_i  \cdot \mathbf{\Phi}_{e,i}(:,N)  \\
    \end{bmatrix} =- \sum_{i=1}^{m} p_{rdc,i} A_i  \begin{bmatrix}  \cos{\psi_i} \\
    {\hat{\boldsymbol{c}}_3}_i  \cdot \mathbf{\Phi}_{e,i}(:,1) \\
    \vdots  \\
    {\hat{\boldsymbol{c}}_3}_i  \cdot \mathbf{\Phi}_{e,i}(:,N)  \\
    \end{bmatrix}  \dot{\boldsymbol{r}}_{dma3,i} \\
    &= - \underbrace{ \sum_{i=1}^{m} p_{rdc,i} A_i  \begin{bmatrix}  \cos{\psi_i} \\
    {\hat{\boldsymbol{c}}_3}_i  \cdot \mathbf{\Phi}_{e,i}(:,1) \\
    \vdots  \\
    {\hat{\boldsymbol{c}}_3}_i  \cdot \mathbf{\Phi}_{e,i}(:,N)  \\
    \end{bmatrix} \begin{bmatrix} C_{sa}(3,3) & [C_{se} \boldsymbol{\Phi}_e(\phi_i)](3,:)
    \end{bmatrix}}_{\boldsymbol{D}_r} \dot{\vec{\mathbf{x}}}     \label{eq:81_2}
\end{align}

Accordingly, the radiation damping matrix, $\boldsymbol{C}$, is expressed as: \begin{equation}
    \boldsymbol{D}_r = \sum_{i=1}^{m} p_{rdc,i} A_i  \begin{bmatrix}  \cos{\psi_i} \\
    {\hat{\boldsymbol{c}}_3}_i  \cdot \mathbf{\Phi}_{e,i}(:,1) \\ 
    \vdots  \\
    {\hat{\boldsymbol{c}}_3}_i  \cdot \mathbf{\Phi}_{e,i}(:,N)  \\
    \end{bmatrix} \begin{bmatrix}1 & [C_{se} \boldsymbol{\Phi}_e(\phi_i)](3,:)
    \end{bmatrix}
\end{equation} where the proposed $\boldsymbol{D}_r$ is a full, non symmetric matrix.

\subsubsection{Derivation of the Generalized Added Inertia Force $\boldsymbol{Q}^{\infty}$:}\label{sec:inertia}
The derivation of the generalized added inertia will yield an expression for the generalized added mass matrix. The added mass pressure is calculated using the complex part of the radiation pressure coefficient as shown in Eq.~\eqref{eq:37_2} \cite{zou2020modeling}: \begin{equation} \label{eq:37_2}
    p_{amc,i} = \frac{\Im \left( p_{rd,i} \right)}{\omega}
\end{equation} such that the added inertia for the VSB WEC is calculated for each mesh element as follows: \begin{equation}
  {F_{\infty,i}} =  p_{amc,i} A_i  
\end{equation} And its corresponding normal vector to the mesh element is expressed as:\begin{align} 
\vec{f}_i^{\infty} &= - {F_{\infty,i}}  \hat{\boldsymbol{c}}_{3,i} \label{eq:36} \\
&= - p_{amc,i} A_i    {\hat{\boldsymbol{c}}}_{3,i} {}^a\ddot{\vec{\boldsymbol{r}}}_{dma3,i}  = - p_{amc,i} A_i    {\hat{\boldsymbol{c}}}_{3,i} \begin{bmatrix}
C_{sa}(3,3) & [C_{se} \boldsymbol{\Phi}_e(\phi_i)](3,:)
\end{bmatrix} \ddot{\vec{\mathbf{x}}} \label{eq:36_2}
\end{align} 

From Eq.~\eqref{eq:32} one can write ${}^a{\ddot{\boldsymbol{r}}_{dma}}$ as: \begin{align}
   {}^a\ddot{\vec{\boldsymbol{r}}}_{dma,i} &= \begin{bmatrix}
C_{sa} & C_{se} \boldsymbol{\Phi}_e(\phi_i)
\end{bmatrix} \ddot{\vec{\mathbf{x}}} \label{eq:37} \\
   {}^a\ddot{\vec{\boldsymbol{r}}}_{dma3,i} &= \begin{bmatrix}
C_{sa}(3,3) & [C_{se} \boldsymbol{\Phi}_e(\phi_i)](3,:)
\end{bmatrix} \ddot{\vec{\mathbf{x}}} 
\end{align} 

The generalized added inertia is the calculated by applying Eq.~\eqref{eq:60_1g} to get: \begin{equation} \label{eq:38}
\boldsymbol{Q}_j^{\infty} = \vec{\boldsymbol{f}}^{\infty} \cdot \frac{\partial \vec{r}_{dma}}{\partial q_j} , \, j = 1,2, ..., 1+N
\end{equation} 

We can then write the added mass generalized force using Eqs. \eqref{eq:36}, \eqref{eq:36_2}, and \eqref{eq:38} as follows: \begin{align}
\boldsymbol{Q}_{1,i}^{\infty} = \vec{\boldsymbol{f}}^{\infty}_i  \cdot \frac{\partial  \vec{\boldsymbol{r}}_{dma} }{\partial \vec{\boldsymbol{r}}_{sa,3}} &= - p_{amc,i} A_i  {}^a\ddot{\vec{\boldsymbol{r}}}_{dma3,i}  {\hat{\boldsymbol{c}}}_{3,i}  \cdot \boldsymbol{1}_3 = - p_{amc,i} A_i  \boldsymbol{1}_3^T {}^a\ddot{\vec{\boldsymbol{r}}}_{dma3,i}  {\hat{\boldsymbol{c}}}_{3,i} \\
 \boldsymbol{Q}_{2,i}^{\infty} = \vec{\boldsymbol{f}}^{\infty}_i \cdot \frac{\partial  \vec{\boldsymbol{r}}_{dma} }{\partial {\eta_1} } &= -p_{amc,i} A_i   {\hat{\boldsymbol{c}}}_{3,i}  {}^a\ddot{\vec{\boldsymbol{r}}}_{dma3,i}  \cdot \mathbf{\Phi}_{e,i}(:,1)  \\
\boldsymbol{Q}_{1+N,i}^{\infty} =  \vec{\boldsymbol{f}}^{\infty}_i \cdot \frac{\partial  \vec{\boldsymbol{r}}_{sa} }{\partial {\eta_N} } &=  -p_{amc,i} A_i   {}^a\ddot{\vec{\boldsymbol{r}}}_{dma3,i}  {\hat{\boldsymbol{c}}}_{3,i} \cdot \mathbf{\Phi}_{e,i}(:,N) 
\end{align} 

Therefore, \begin{align}
    \mathbit{Q}^{\infty} &= -\sum_{i=1}^{m} p_{amc,i} A_i \begin{bmatrix}     \cos{\psi_i}  \\
      {\hat{\boldsymbol{c}}_3}_i  \cdot \mathbf{\Phi}_{e,i}(:,1) \\ 
    \vdots  \\
    {{}^a\ddot{\vec{\boldsymbol{r}}}_{dma3,i} \hat{\boldsymbol{c}}_3}_i  \cdot \mathbf{\Phi}_{e,i}(:,N)  \\
    \end{bmatrix}  {}^a\ddot{\vec{\boldsymbol{r}}}_{dma3,i}  \\
    &= - \underbrace{ \sum_{i=1}^{m} p_{amc,i} A_i \begin{bmatrix}     \cos{\psi_i}  \\
      {\hat{\boldsymbol{c}}_3}_i  \cdot \mathbf{\Phi}_{e,i}(:,1) \\ 
    \vdots  \\
    {\hat{\boldsymbol{c}}_3}_i  \cdot \mathbf{\Phi}_{e,i}(:,N)  \\
    \end{bmatrix}  \begin{bmatrix} C_{sa}(3,3) & [C_{sa} \boldsymbol{\Phi}_e(\phi_i)](3,:)
\end{bmatrix}}_{\boldsymbol{M}_{\infty} } \ddot{\vec{\mathbf{x}}}   \label{eq:44_2}
\end{align}

From Eq.~\eqref{eq:44_2} we can conclude that the generalized added mass for heave-only motion equals:
\begin{align}
    \boldsymbol{M}_{\infty} &= \sum_{i=1}^{m} p_{amc,i} A_i \begin{bmatrix}    \cos{\psi_i}  \\
      {\hat{\boldsymbol{c}}_3}_i  \cdot \mathbf{\Phi}_{e,i}(:,1) \\ 
    \vdots  \\
    {\hat{\boldsymbol{c}}_3}_i  \cdot \mathbf{\Phi}_{e,i}(:,N)  \\
    \end{bmatrix}  \begin{bmatrix}
1 & [C_{se} \boldsymbol{\Phi}_e(\phi_i)](3,:)
\end{bmatrix}    \label{eq:45}
\end{align}

The proposed generalized added mass derived in Eq.~\eqref{eq:45} is a full non-symmetric matrix, i.e, the equation of motion becomes dynamically coupled when adding the added mass to the system mass matrix $\boldsymbol{M}$.

Finally, the proposed equation of motion for regular waves can be expressed as:\begin{equation} \label{eq:105_3}
\tcbhighmath{ \underbrace{(\boldsymbol M + \boldsymbol M_{\infty})}_{\Tilde{\boldsymbol{M}}(t)} \ddot{\mathbf{x}} + \underbrace{(\boldsymbol{D} + \boldsymbol{D}_r)}_{\Tilde{\boldsymbol{D}}(t)} \dot{\mathbf{x}} +  \underbrace{\left( \boldsymbol{K} +  \boldsymbol{K}_h \right)}_{\Tilde{\boldsymbol{K}}(t)} {\mathbf{x}}  =  \boldsymbol Q^c (t)
+ \boldsymbol Q^{ext}(t) + \boldsymbol Q^{pto}(t)} 
\end{equation} which can be written in the following form: \begin{equation}
    \Tilde{\boldsymbol{M}}(t) \ddot{\mathbf{x}} + \Tilde{\boldsymbol{D}}(t) \dot{\mathbf{x}} +  \Tilde{\boldsymbol{K}}(t) {\mathbf{x}}  =  \boldsymbol Q^c (t)
+ \boldsymbol Q^{ext} (t)+ \boldsymbol Q^{pto}(t)
\label{eq:EOM} 
\end{equation} Noting that the expression of the constraint force $\boldsymbol Q^c$
on VSB WEC shell is derived and demonstrated in Ref. \cite{shabara2022Dynamics}. 

\subsection{Generalized Hydrodynamic Forces For Irregular Waves}
An irregular wave is considered as a set of superimposed regular waves with phase shifts "$\hat{\phi}$", and these regular waves have different frequencies and amplitudes. The amplitude of each of these regular waves is calculated using the Wave Spectrum (short for Wave Energy Density Spectrum). The Wave Spectrum is a function that specifies the amount of energy contained in each wave frequency \cite{wacher2010mathematical}. In this paper, the Bretchsneider Spectrum \cite{brodtkorb2000wafo} shown in Eq.~\eqref{eq:65} is used \begin{equation} \label{eq:65}
    S(f) = \frac{5}{16} \frac{f_p^4}{f^5}H_s^2e^{{-\frac{5}{4}\frac{f_p^4}{f^4}}}
\end{equation}

where $H_s$ is the significant wave height and $f_p = 1/T_p= \omega_p/(2 \pi)$ is the frequency of a particular wave component.

The simulated wave spectrum consisted of an array of $n$ equidistant frequencies with step $\Delta f$, such that the $i^{th}$ frequency is calculated using Eq. (\ref{eq:66}) \cite{wacher2010mathematical}
\begin{equation} \label{eq:66}
    f_i =  f_1 + (i-1) \Delta f  
\end{equation}
The wave height $\eta$ is then calculated as $\eta = \sqrt{2S(f)\Delta f}$  \cite{wacher2010mathematical}.

\subsubsection{Derivation of the $\boldsymbol{Q}^{hst}$ and $\boldsymbol{Q}^{\infty}$ Forces for Irregular Waves:}
The derivation of the expressions for the generalized hydrostatic force $\boldsymbol{Q}^{hst}$ and stiffness $\boldsymbol{K}_h$ does not change from regular to irregular waves, i.e., the exact derivation in subsection \ref{sec:buoyant} applies.

Similarly, the derivation of the expressions for the generalized inertia $\boldsymbol{Q}^{\infty}$ and the added mass matrix $\boldsymbol{M}_\infty$ presented in section \ref{sec:inertia} does not change from regular to irregular waves, except that for the irregular waves, the calculations are based on the largest wave frequency present in the wave spectrum ($\omega_{\text{max}}$). \\ 

%In the next subsections the excitation and radiation hydostatic forces for irregular waves are derived. 

% Since the buoy is constrained from any rotations, the $\mathfrak{A}$ and $\mathfrak{B}$ reference frames axes are parallel to each other and $\hat{a}_3 \equiv \hat{b}_3 $.   
\subsubsection{Derivation of the $\boldsymbol{Q}^{ext}$ Force:}

The time domain excitation force on the $i^{th}$ panel on the buoy surface can be expressed as \cite{zou2020modeling}: \begin{equation} \label{eq:54}
 p_{ext,i} = \sum_{j=1}^{N_w} \Re \left( p_{ex,i}(\omega_j) \eta(\omega_j) e^{\hat{i} \left( \omega_j t + \hat{\phi}_j \right)} \right)
\end{equation} where $N_w$ is the number of superimposed frequencies in the wave, $\hat{\phi}$ is a random phase shift, and $\hat{i}=\sqrt{-1}$, following the same procedure as for the regular wave derivation in subsection \ref{sssec:exh} but with replacing Eq.~\eqref{eq:29_1} with Eq.~\eqref{eq:54} in the derivation. This yields a similar expression for the generalized hydro force related to the excitation and the hydrostatic forces "$\boldsymbol{Q}^{ext}$" as in Eq.~\eqref{eq:81}. Accordingly, the excitation force coefficient for VSB WEC for the $j^{th}$ frequency is expressed as: \begin{equation}
    \boldsymbol{Ex}(\omega_j,t) =  \sum_{i=1}^m p_{ex,i}(\omega_j) A_i  \begin{bmatrix}  \cos{\psi_i} \\
    {\hat{\boldsymbol{c}}_3}_i  \cdot \mathbf{\Phi}_{e,i}(:,1) \\ 
    \vdots  \\
    {\hat{\boldsymbol{c}}_3}_i  \cdot \mathbf{\Phi}_{e,i}(:,N)  \\
    \end{bmatrix} \label{eq:89_3}
\end{equation} Note that for the two-way FSI simulations, "$\boldsymbol{Ex}(\omega,t)$" is time and frequency dependent as it's a function of the buoy's shape. Finally, the generalized excitation force for VSB WECs in irregular waves is calculated as: \begin{equation}
    \boldsymbol{Q}^{ext}(t) = \sum_{j=1}^{N_\omega} \Re \left( \boldsymbol{Ex}_j \; \eta_j \; e^{\hat{i} (\omega_j t + \hat{\phi}_j )} \right) \end{equation}

\subsubsection{Derivation of the $\boldsymbol{Q}^{rad}$ Force:}
The radiation damping force on the $i^{th}$ mesh element acting on the normal direction is computed as: \begin{align} 
\vec{f}_i^{rad} &= - {F_{rad,i}}  \hat{\boldsymbol{c}}_{3,i} \label{eq:88_4} \\
    F_{rad,i} &= p_{rdt,i} A_i \label{eq:100_1}
\end{align} where $p_{rdt,i}$ is the time domain radiation pressure, and it is derived in Ref. \cite{zou2020modeling} as: \begin{align}
    p_{rdt.i} &= \frac{2}{\pi} \int_0^t \sum_{\omega_{min}}^{\omega_{max}} \left( p_{rdc,i} \cos{\omega (t-\tau)} \: \Delta \omega {}^a\dot{\vec{\boldsymbol{r}}}_{dma3,i}(\tau) \right) d  \tau \\
    & \approx \frac{2}{\pi} \sum_0^t \sum_{\omega_{min}}^{\omega_{max}} \left( p_{rdc,i} \cos{\omega (t-\tau)} \: \Delta \omega {}^a\dot{\vec{\boldsymbol{r}}}_{dma3,i}(\tau) \right) \Delta  \tau \label{eq:101_1}
\end{align} From Eq.~\eqref{eq:60_1g} the generalized hydro force can be written in the following form:
\begin{equation} \label{eq:60_2} 
  \mathbit{Q}_j^{rad} = \vec{\boldsymbol{f}}^{rad} \cdot \frac{\partial \vec{r}_{dma}}{\partial q_j} , \, j = 1,2, ..., 1+N
\end{equation} And accordingly we can then write the generalized radiation force using Eqs.~\eqref{eq:60_2} and~\eqref{eq:88_4} as follows: \begin{align}
\mathbit{Q}_{1,i}^{rad} &= \vec{\boldsymbol{f}}^{rad}_i  \cdot \frac{\partial  \vec{\boldsymbol{r}}_{dma} }{\partial \vec{\boldsymbol{r}}_{sa,3}} = - {F_{rad,i}}  \hat{\boldsymbol{c}}_{3,i}   \cdot \boldsymbol{1}_3 = - {F_{rad,i}} \boldsymbol{1}_3^T     \hat{\boldsymbol{c}}_{3,i}  \\
 \mathbit{Q}_{2,i}^{rad} &= \vec{\boldsymbol{f}}^{rad}_i \cdot \frac{\partial  \vec{\boldsymbol{r}}_{dma} }{\partial {\eta_1} } =  - {F_{rad,i}}  \hat{\boldsymbol{c}}_{3,i}   \cdot \mathbf{\Phi}_{e,i}(:,1)  \\
\mathbit{Q}_{1+N,i}^{rad} &=  \vec{\boldsymbol{f}}^{rad}_i \cdot \frac{\partial  \vec{\boldsymbol{r}}_{sa} }{\partial {\eta_N} } =   - {F_{rad,i}}  \hat{\boldsymbol{c}}_{3,i}  \cdot \mathbf{\Phi}_{e,i}(:,N) 
\end{align} Finally, the generalized radiation force is expressed as: \begin{equation}
     \mathbit{Q}^{rad} = -\sum_{i=1}^{m} {F_{rad,i}}   \begin{bmatrix}     \cos{\psi_i}  \\
      {\hat{\boldsymbol{c}}_3}_i  \cdot \mathbf{\Phi}_{e,i}(:,1) \\ 
    \vdots  \\
  {\hat{\boldsymbol{c}}_3}_i  \cdot \mathbf{\Phi}_{e,i}(:,N)  \\
    \end{bmatrix}  \label{eq:91_1}
\end{equation}

When substituting Eqs~\eqref{eq:101_1} and \eqref{eq:100_1} in the expression above, one can notice that the generalized radiation damping force is computed using the convolution of the multiplication of the retardation function and the velocity states as shown in Eq.~\eqref{eq:55_3} \begin{align}
    \boldsymbol{Q}^{rad} &= - \sum_{i=1}^m \int_0^t  \left( \boldsymbol{K}_i(t - \tau) \dot{\mathbf{x}}(\tau) \right) d\tau \label{eq:55_3} \\
    & \approx -  \sum_{i=1}^m \sum_0^t  \left( \boldsymbol{K}_i (t - \tau) \dot{\mathbf{x}}(\tau) \right) \Delta \tau \label{eq:55_2}
\end{align} where the retardation function (memory function) $\boldsymbol{K} \in \mathbb{R}^{(N+1) \times (N+1)}$  for the $i^{th}$ mesh element for heave-only motion can be calculated as: \begin{align}
    \boldsymbol{K}_i(t) &= \frac{2}{\pi} \int_0^\infty \boldsymbol{D}_i(\omega) \cos{(\omega t)} d\omega = \frac{2}{\pi} \sum_{\omega_{min}}^{\omega_{max}} \boldsymbol{D}_i(\omega) \cos{(\omega t)} \Delta \omega
\end{align}
Accordingly, the generalized radiation damping coefficient $\boldsymbol{D}_i(\omega) \in \mathbb{R}^{(N+1) \times (N+1)}$ for the $i^{th}$ mesh element is  expressed as: \begin{equation}
    \boldsymbol{D}_i(\omega) = p_{rdc,i}(\omega) A_i  \begin{bmatrix}     \cos{\psi_i}  \\
      {\hat{\boldsymbol{c}}_3}_i  \cdot \mathbf{\Phi}_{e,i}(:,1) \\ 
    \vdots  \\
  {\hat{\boldsymbol{c}}_3}_i  \cdot \mathbf{\Phi}_{e,i}(:,N)  \\ \end{bmatrix}  \begin{bmatrix}
1 & [C_{se} \boldsymbol{\Phi}_e(\phi_i)](3,:)
\end{bmatrix} \label{eq:100_3}
\end{equation} 

To calculate the generalized radiation force apply Eq.~\eqref{eq:91_1}, which requires a convolution calculation that can be computationally expensive. Another way is to compute the generalized radiation damping coefficients vector $\boldsymbol{D}$ for each frequency, then interpolate the result arrays with respect to the particular wave frequency using the MATLAB function \textit{interp1} to get an averaged generalized radiation damping coefficient ${\boldsymbol{D}_r}(t)$. The generalized radiation damping force is then calculated as  $\boldsymbol{Q}_{rad}={\boldsymbol{D}_r}(t) \dot{\mathbf{x}}$.

An alternative method is to use the state space representation similar to the model described in Eq.~\eqref{eqXr}, noting that each state in the $\mathbf{x}$ vector will have its corresponding linear time variant (LTV) radiation state space model. Finally, Eq.~\eqref{eq:91_1} is used in this work as the convolution provides the more accurate solution.\\

Finally, the proposed equation of motion for VSB WEC in irregular waves can be expressed as: \begin{equation} \label{eq:105_4}
\tcbhighmath{ \underbrace{(\boldsymbol M + \boldsymbol M_{\infty})}_{\Tilde{\boldsymbol{M}}(t)} \ddot{\mathbf{x}} + \boldsymbol{D}\dot{\mathbf{x}} +  \underbrace{\left( \boldsymbol{K} +  \boldsymbol{K}_h \right)}_{\Tilde{\boldsymbol{K}}(t)} {\mathbf{x}}  =  \boldsymbol Q^c (t)
+ \boldsymbol Q^{ext}(t)+ \boldsymbol Q^{rad}(t) + \boldsymbol Q^{pto}(t)}
\end{equation} which can be written in the following form: \begin{equation}
    \Tilde{\boldsymbol{M}}(t) \ddot{\mathbf{x}} + {\boldsymbol{D}} \dot{\mathbf{x}} +  \Tilde{\boldsymbol{K}}(t) {\mathbf{x}}  =  \boldsymbol Q^c (t)
+ \boldsymbol Q^{ext} (t)+  \boldsymbol Q^{rad} (t) +\boldsymbol Q^{pto} (t)
\label{eq:EOM_2w} 
\end{equation}

\section{The Equation of Motion of VSB WECs:} \label{sec:FSI}
\begin{figure}[h]
\centering
\includegraphics[scale=0.9]{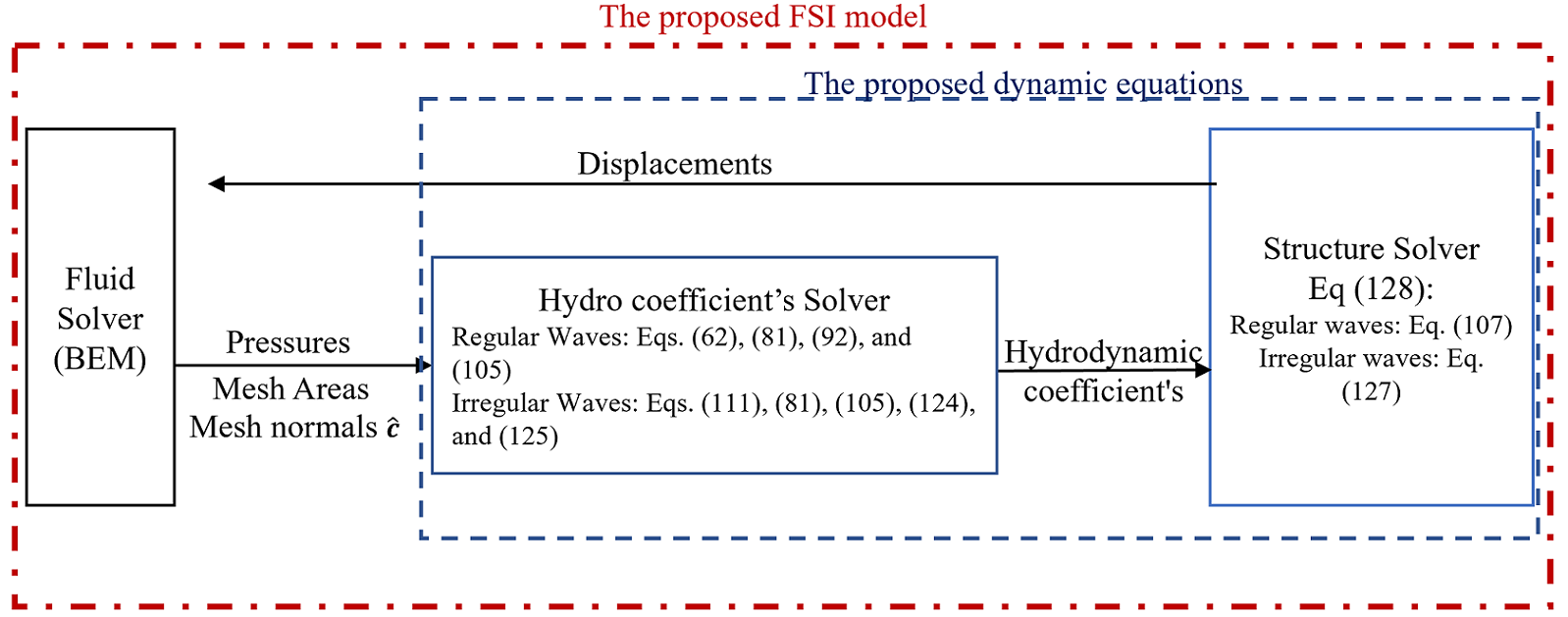}
\caption{The proposed FSI model for VSB WECs}
\label{fig:FSI}
\end{figure}

The proposed equation of motion for the VSB WEC in this work can be expressed in a more general form similar to Eq.~\eqref{eq:Cummins} as:
\begin{equation} 
{\boldsymbol{M}}  \ddot{\mathbf{x}}(t) {+} \boldsymbol{D} \dot{\mathbf{x}}(t) {+} \boldsymbol{K} {\mathbf{x}}(t) = \overbrace{\int_{-\infty}^{\infty}\boldsymbol{H}_{ext}(\tau)\eta(t {-} \tau,\mathbf{x})d\tau}^\text{Generalized excitation force $\boldsymbol{Q}^{{ext}}$}  + \boldsymbol{Q}^{c} +\boldsymbol{Q}^{{hst}}(t) \overbrace{-{\boldsymbol{M}}_\infty(t)  \ddot{\mathbf{x}}(t) - \int_{-\infty}^{t} \boldsymbol{K}(\tau)\dot{\mathbf{x}}(t - \tau)d\tau}^\text{Generalized radiation force $\boldsymbol{Q}^{'rad}$} - \boldsymbol{Q}^{pto}   \label{EoMF_VSB}
\end{equation} where $\boldsymbol{H}_{ext}$ is the impulse response function defining the generalized excitation force for heave. It is possible to show that this equation reduces to Eq.~\eqref{eq:105_3} for regular waves, and to Eq.~\eqref{eq:105_4} for irregular waves. Fig. \ref{fig:FSI} shows a schematic of the proposed FSI model in which the BEM code solves for the hydrodynamic and hydrostatic pressures. Then the pressures and mesh details are passed to a function that calculated the hydodynamic coefficients based on the proposed equations in section \ref{sec:const_forced}. Then the structure solver takes these hydrodynamic coefficients as an input and calculates their corresponding generalized forces based on the proposed expressions in section \ref{sec:const_forced}; the structure solver then calculates the system's response using Eq. \eqref{EoMF_VSB} to get $\mathbf{x}$, and $\dot{\mathbf{x}}$. This process can be done in a two-way, or a one-way FSI scheme. 

Next section proposes a rigorous method to calculate averaged hydro coefficients for the one-way FSI model based on Reynolds averaging method. 

\section{One Way FSI Equation of Motion Derivation} \label{sec:one-two_way}
The two-way FSI simulations are computationally expensive even when running simulations for regular waves. A one-way FSI model that can replace the two-way FSI simulations without compromising the accuracy of the solution is highly desirable. To get the required averaged quantities, we propose using the Reynolds averaging method, this method is used in CFD turbulence modeling to obtain the Reynolds-averaged Navier–Stokes equations (RANS equations) \cite{shabara2018numerical,goccmen2016wind,durbin2021basics}. 

The time-dependent coefficients $\Tilde{\boldsymbol{M}}$, $\Tilde{\boldsymbol{D}}$, $\Tilde{\boldsymbol{K}}$, in Eqs.~\eqref{eq:EOM} and \eqref{eq:EOM_2w} can be expressed as a summation of an averaged value (statistical mean) and a fluctuating value (from the mean); this is called Reynolds decomposition \cite{durbin2021basics}, such that:
\begin{align}
    \Tilde{\boldsymbol{M}}(t) &= \overline{\boldsymbol{M}} + \boldsymbol{M}'(t) \longrightarrow  \overline{\Tilde{\boldsymbol{M}}} = \overline{\overline{\boldsymbol{M}} + {\boldsymbol{M}'}}  = \overline{\overline{\boldsymbol{M}}} = \overline{\boldsymbol{M}} = \lim_{T \rightarrow \infty} \int^T_0 \Tilde{\boldsymbol{M}}(t) \; dt  \\
    \Tilde{\boldsymbol{D}}(t) &= \overline{\boldsymbol{D}} + {\boldsymbol{D}'}(t)\longrightarrow \overline{\Tilde{\boldsymbol{D}}} =\overline{\boldsymbol{D}} = \lim_{T \rightarrow \infty} \int^T_0 \Tilde{\boldsymbol{D}}(t) \; dt  \\
    \Tilde{\boldsymbol{K}}(t) &= \overline{\boldsymbol{K}} + {\boldsymbol{K}'}(t)\longrightarrow \overline{\Tilde{\boldsymbol{K}}} =\overline{\boldsymbol{K}}= \lim_{T \rightarrow \infty} \int^T_0 \Tilde{\boldsymbol{K}}  (t) \;dt 
    \boldsymbol{}
\end{align} A similar approach can be followed to obtain averaged values for the generalized excitation force coefficients for regular and irregular waves expressed in Eqs.~\eqref{eq:53_3} and \eqref{eq:89_3}. The resulted generalized averaged excitation force for irregular waves can be expressed as:\begin{equation}
 \overline{\boldsymbol{Q}^{ext}}(t) = \sum_{j=1}^{N_\omega} \Re \left( \overline{\boldsymbol{Ex}}(\omega_j) \eta_j e^{\hat{i} \left( \omega_j t + \hat{\phi}_j \right)} \right)
 \end{equation} and the generalized averaged radiation force $\overline{\boldsymbol{Q}^{rad}}$ is calculated as: \begin{equation} \label{eq:99_3}
   \overline{\boldsymbol{Q}^{rad}}(t)  \approx -  \sum_{i=1}^m \sum_0^t  \left( \overline{\boldsymbol{K}_i} (t-\tau) \dot{\mathbf{x}}(\tau) \right) \Delta t   \end{equation} where $\overline{\boldsymbol{K}_i}$ is the average generalized retardation function calculated based on the averaged generalized radiation damping coefficient $\overline{\boldsymbol{B}_i}(\omega)$.

A proper one-way FSI simulation should use a proper average value, and the Equations of Motion \eqref{eq:EOM} and \eqref{eq:EOM_2w} can be expressed as:  \begin{align}
        \overline{\boldsymbol{M}} \ddot{\mathbf{x}} + \overline{\boldsymbol{D}} \dot{\mathbf{x}} +  \overline{\boldsymbol{K}} {\mathbf{x}}  &=  \overline{\boldsymbol Q^c} (t)
+ \overline{\boldsymbol{Q}^{ext}}(t) + \boldsymbol Q^{pto}(t) \label{eq:EOM_1w} \\
        \overline{\boldsymbol{M}} \ddot{\mathbf{x}} + {\boldsymbol{D}} \dot{\mathbf{x}} +  \overline{\boldsymbol{K}} {\mathbf{x}}  &=  \overline{\boldsymbol Q^c} (t)
+\overline{\boldsymbol{Q}^{ext}}(t) +\overline{\boldsymbol{Q}^{rad}}(t) + \boldsymbol Q^{pto}(t) \label{eq:EOM_1w_2} 
\end{align} respectively, where \begin{align}
    \overline{\boldsymbol{M}} &= \boldsymbol{M}+ \overline{{\boldsymbol{M}_\infty}}  \\
    \overline{\boldsymbol{D}} &= \boldsymbol{D}+ \overline{{\boldsymbol{D}_r}} \\ 
    \overline{\boldsymbol{K}} &= \boldsymbol{K}+ \overline{{\boldsymbol{K}_h}} 
\end{align} and generalized averaged constraint force $\overline{\boldsymbol Q^c}$ is calculated based on the generalized averaged non-conservative forces acting on the WEC system \cite{shabara2022Dynamics}.  

To solve the equations above, average values for generalized added mass, hydrodynamic damping, hydrostatic and excitation force coefficient matrices are needed. For both regular and irregular waves, the VSB oscillates about its initial non-deformed shape. i.e., the frequency domain pressure coefficients obtained at the initial time step (from the non-deformed shape) are be used to calculate the averaged generalized values for required averaged quantities, this will be justified in section \ref{sec:results}.

\section{Model Validation} \label{sec:valid}
The validation of generalized structural coefficients in Eq.~\eqref{eq:50_2}, Rayleigh-Ritz approximation in Eq.~\eqref{eq:5}, and the Legendre polynomials for the mode-shapes in Eq.~\eqref{eq:6} are validated in \cite{shabara2022Dynamics}. This section presents the validation for the derived hydrodynamic coefficients and forces for both regular and irregular waves. Also, the two-way FSI time interval decency test is presented.    

\subsection{Hydrodynamic Model Validation}
The result of the derived FSI model is compared to the results obtained using Cummins equation \eqref{eq:Cummins} for regular and irregular waves. Cummins equation treats the WEC as a rigid undeformable body; on the other hand, the developed FSI model does not have the rigidity assumption. To compare the developed model with the result derived by Cummins equation, the material properties of the VSB WEC are changed to simulate a rigid shell buoy. Accordingly, the Poisson's ratio is set to zero "$\nu = 0$", and the young's modulus is set to 10 GPa. Also, the shell radius is set to $r = 2 m$.

\begin{figure}
     \centering
     \begin{subfigure}{0.46\textwidth}
         \centering
         \includegraphics[width=\textwidth]{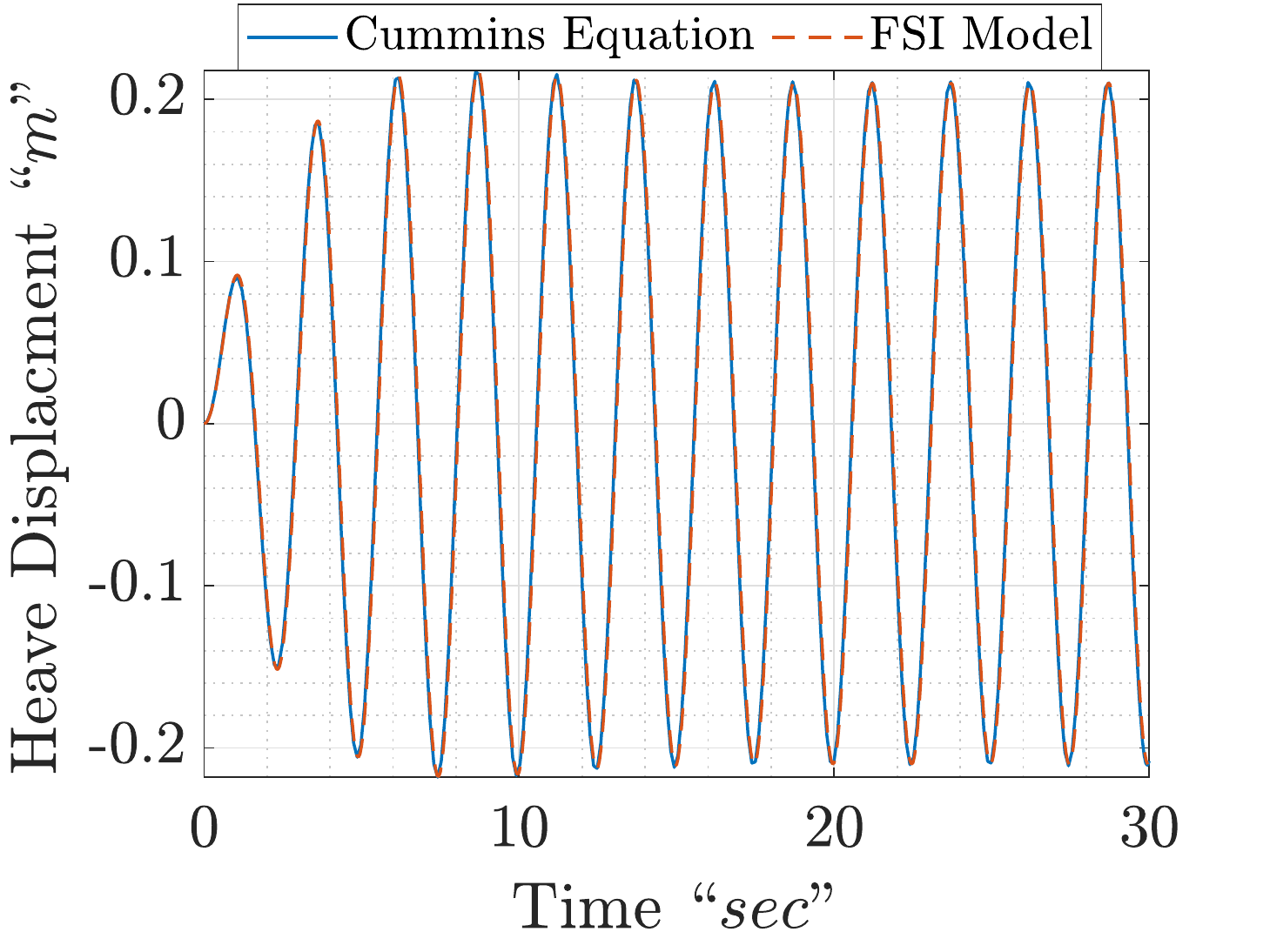}
         \caption{Heave Displacement}
         \label{fig:valid_1a}
     \end{subfigure}
     \hfill
     \begin{subfigure}{0.46\textwidth}
         \centering
         \includegraphics[width=\textwidth]{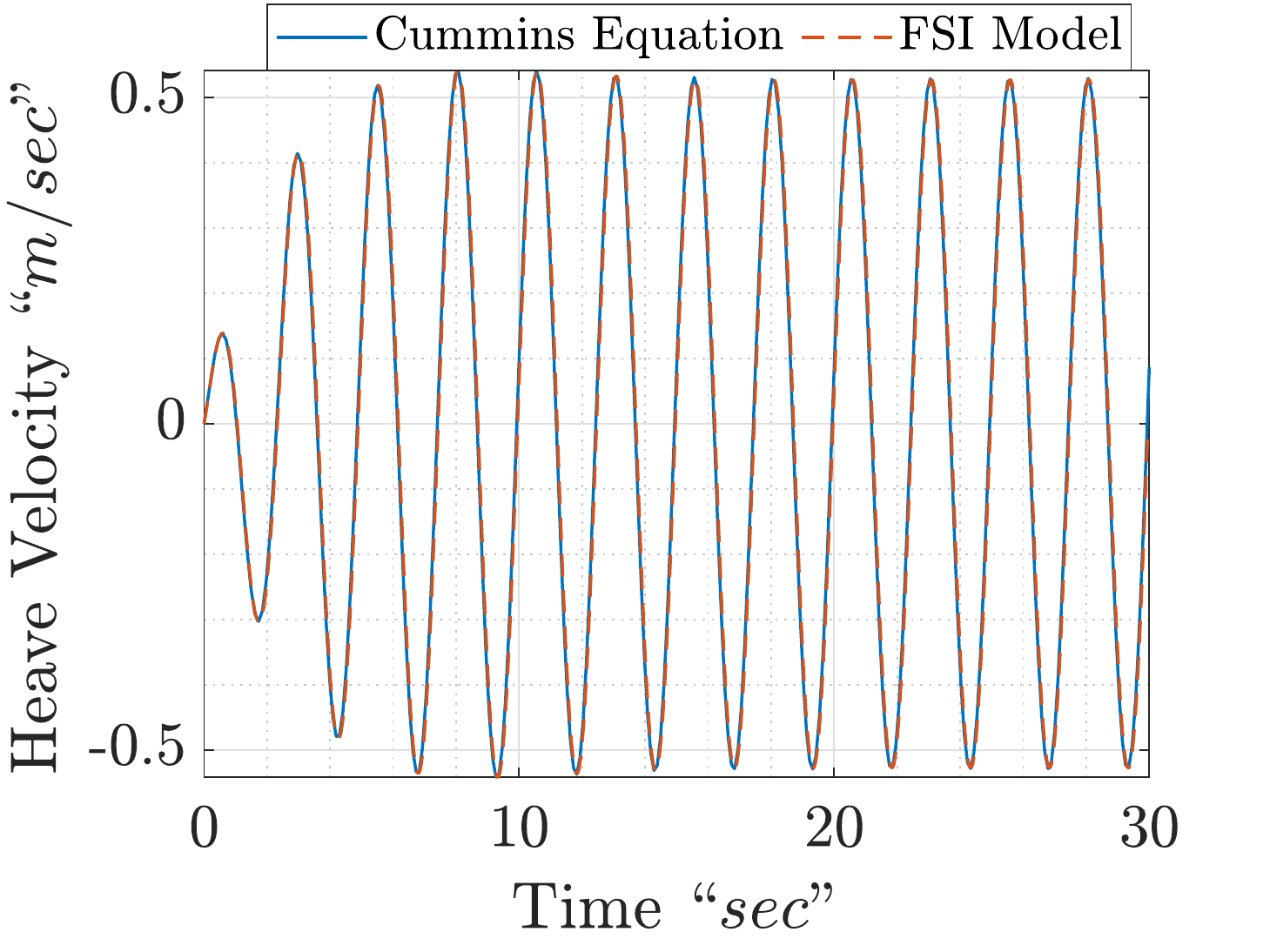}
         \caption{Velocity Displacement}
         \label{fig:valid_1b}
     \end{subfigure}
    \caption{Trajectories resulted from Cummins equation and the two-way FSI model for regular waves}%
    \label{fig:valid_reg}%
\end{figure}

Fig.~\ref{fig:valid_reg} shows the trajectories resulting from both the derived equation of motion vs Cummins equation in a generic regular wave. It can be seen that the trajectories are almost identical and resulted in the harvested energy by the FSI model to produced slightly less energy compared to Cummins Equation by a factor of 0.5\%.  

\begin{figure}
     \centering
     \begin{subfigure}{0.46\textwidth}
         \centering
         \includegraphics[width=\textwidth]{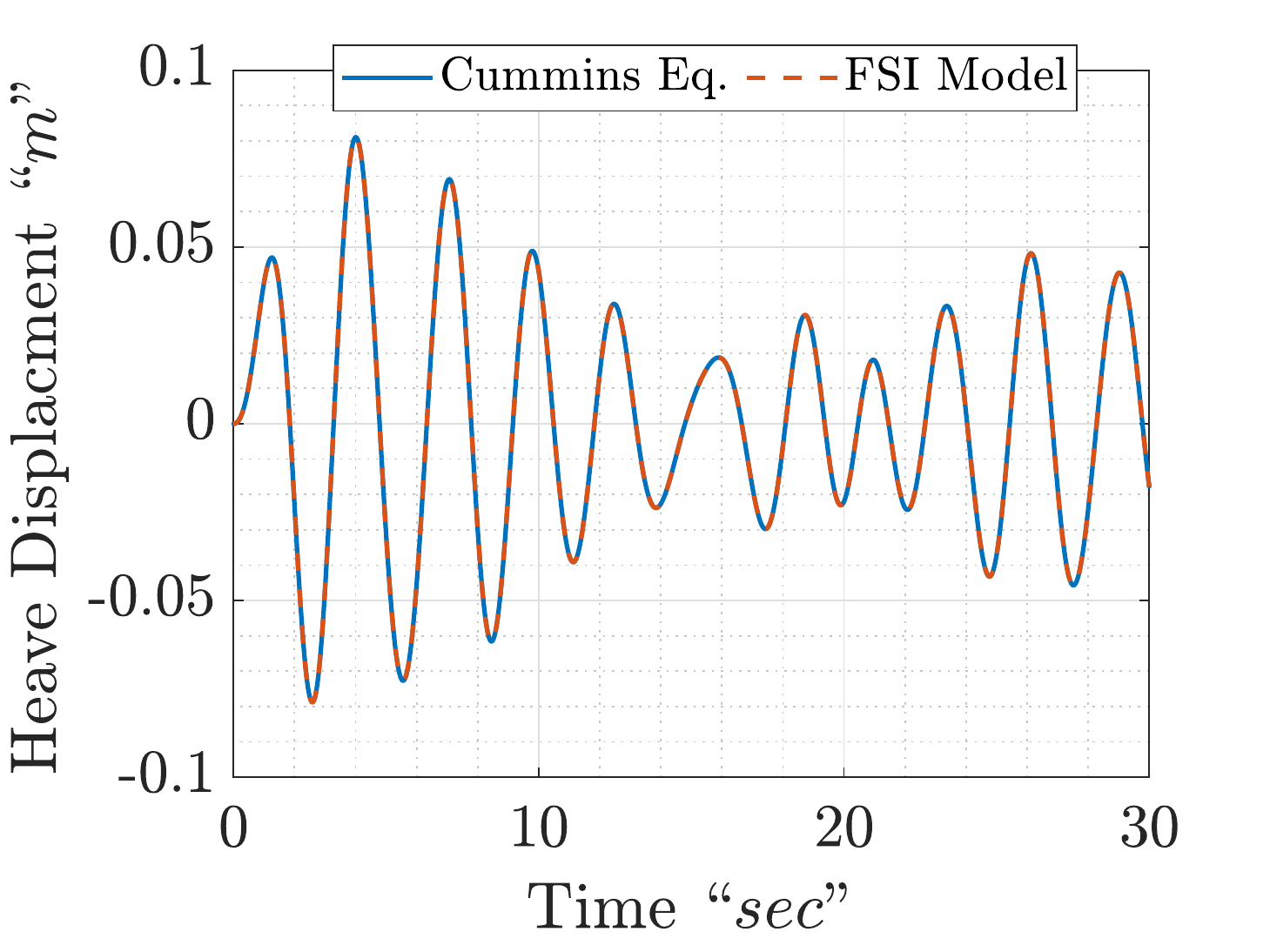}
         \caption{Heave Displacement}
         \label{fig:valid_1c}
     \end{subfigure}
     \hfill
     \begin{subfigure}{0.46\textwidth}
         \centering
         \includegraphics[width=\textwidth]{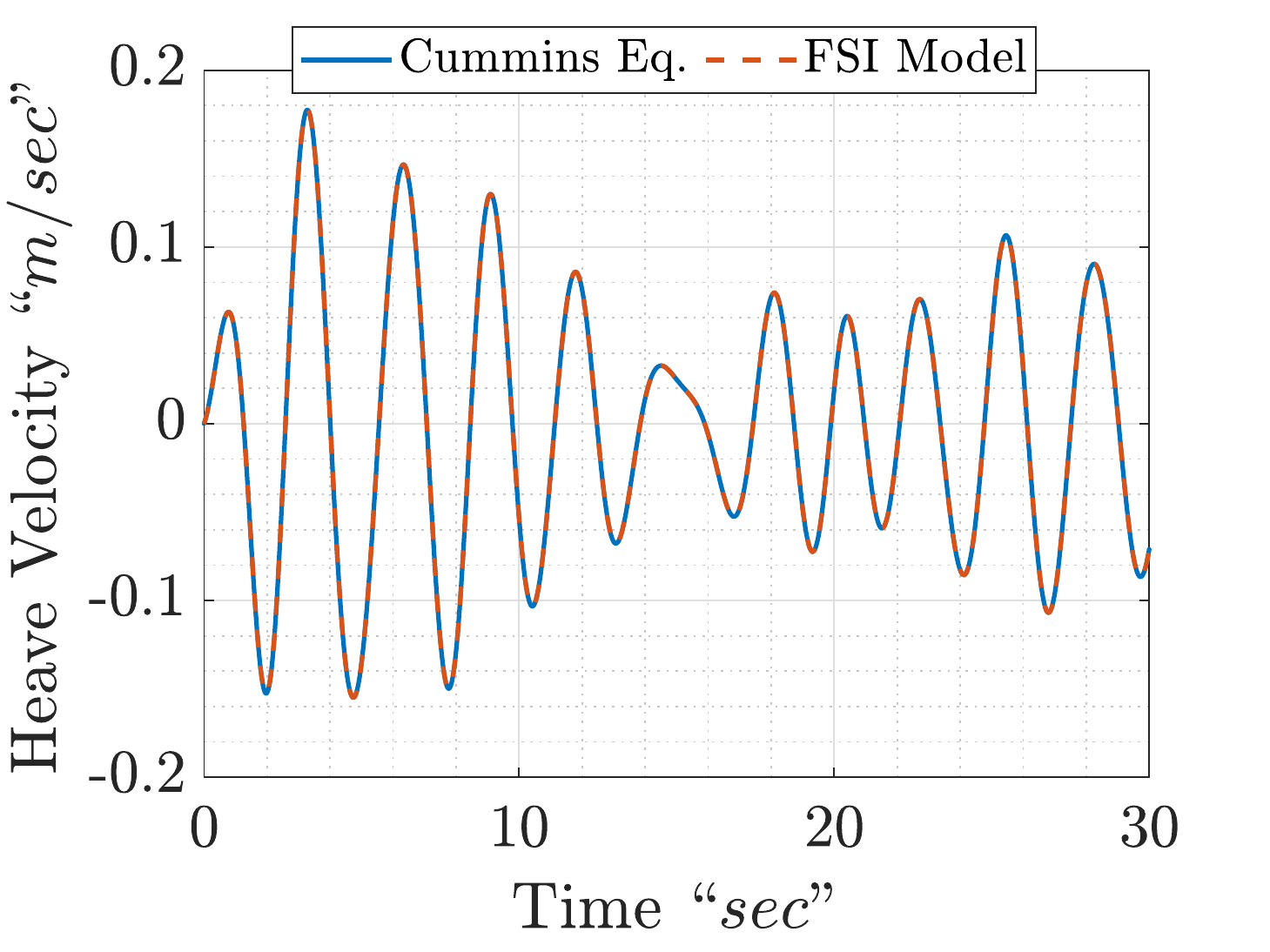}
         \caption{Velocity Displacement}
         \label{fig:valid_1d}
     \end{subfigure}
    \caption{Responses resulted From Cummins equation and the two-way FSI model for irregular waves}%
    \label{fig:valid_irreg}%
\end{figure}

Similarly, Fig.\ref{fig:valid_irreg} shows the waveform obtained from Cummins equation and the developed two-way FSI model in irregular waves. The calculation of the radiation force in both models was done using the convolution technique. It can be noticed that the error between the displacement and velocity responses is negligible. Also, the error in the harvested energy was bounded by $\pm0.5 \%$.

\subsection{FSI intervals Validation}
A multiple-step time integration method is used, namely the six-stage, fifth-order, Runge-Kutta method. To reduce the computation time, the developed FSI model calls NEMOH every preset time internal to update the hydrodynamic/static pressures based on the instantaneous shape of the buoy; multiple time interval sizes were tested ranging from 0.05 to 0.3 seconds. Table~\ref{tab:1} shows the peak-to-peak heave displacement results for various time intervals. It can be seen that the change in the pk-pk displacement for time steps $0.1$ and $0.05$ is negligible. Accordingly, the time step 0.05 seconds was used in the two-way FSI simulations. 
\begin{table}
\caption{Dependency of heave displacement on time intervals for regular waves}\label{tab:1}
\begin{tabular}{|c|c|} 
\hline
\textbf{Time Interval} & \textbf{displacement (m pk-pk)} \\ \hline
0.05                   & 0.4936                          \\ \hline
0.1                    & 0.4915                          \\ \hline
0.2                    & 0.4873                          \\ \hline
0.3                    & 0.4832                          \\ \hline
\end{tabular}
\end{table}

It is worth noting that using a large FSI time interval (eg: 0.3 seconds) for irregular waves would yield a divergence in the convolution integral of the generalized radiation damping force "Eq.~\eqref{eq:99_3}". 

\section{Results and Discussion} \label{sec:results} 
The tested sea conditions and simulation times use the recommended values listed in Table \ref{tab:2}. These values were recommended in Ref.~\cite{abdelkhalik2016control}. Both regular and irregular waves were tested; the irregular waves have a frequency range between 0.1 and 6 rad/s with a number of frequencies  $N_\omega =$ 256. The power take-off unit in this work uses a passive loading control with a damping coefficient $c=8000$ Ns/m, where the damping coefficient is the ratio between the excitation force and the heave velocity of the WEC.

% \usepackage{multirow}
\begin{table}[]
\caption{Tested Wave Conditions \cite{abdelkhalik2016control}}\label{tab:2}
\begin{tabular}{|c|c|c|c|c|}
\hline
\multirow{2}{*}{\textbf{ID}} & \multirow{2}{*}{$T_P$ (s)} & \multirow{2}{*}{$H_s$ (m)} & \multirow{2}{*}{Duration (s)}  \\
 &    &           &                              \\ \hline
RS06 & 2.5 & 0.194  & 300                    \\ \hline
RS07 & 3   & 0.278  & 300                  \\ \hline
RS08 & 3.5 & 0.37   & 300                   \\ \hline
RS09 & 4   & 0.464  & 300                    \\ \hline
RS10 & 4.5 & 0.556  & 300                   \\ \hline
RS11 & 5   & 0.646  & 300                  \\ \hline
RS12 & 6   & 0.8222 & 300                    \\ \hline
RS13 & 7   & 0.992  & 360                  \\ \hline
RS14 & 8   & 1.158  & 360                  \\ \hline
\end{tabular}
\end{table}

This section is divided into three subsections, subsection \ref{subsec:reg} discusses the performance of the VSB WEC in regular waves in two-way FSI environment, and the second subsection is a comparison between the one-way and two-way FSI results for VSB WEC in regular waves where it was proven that the one-way FSI model can replace the two-way FSI model without compromising the results' accuracy. The third subsection discusses the simulations of the VSB WEC in irregular waves.

\subsection{Regular Waves} \label{subsec:reg}
\begin{figure}[h]
\centering
\includegraphics[scale=0.6]{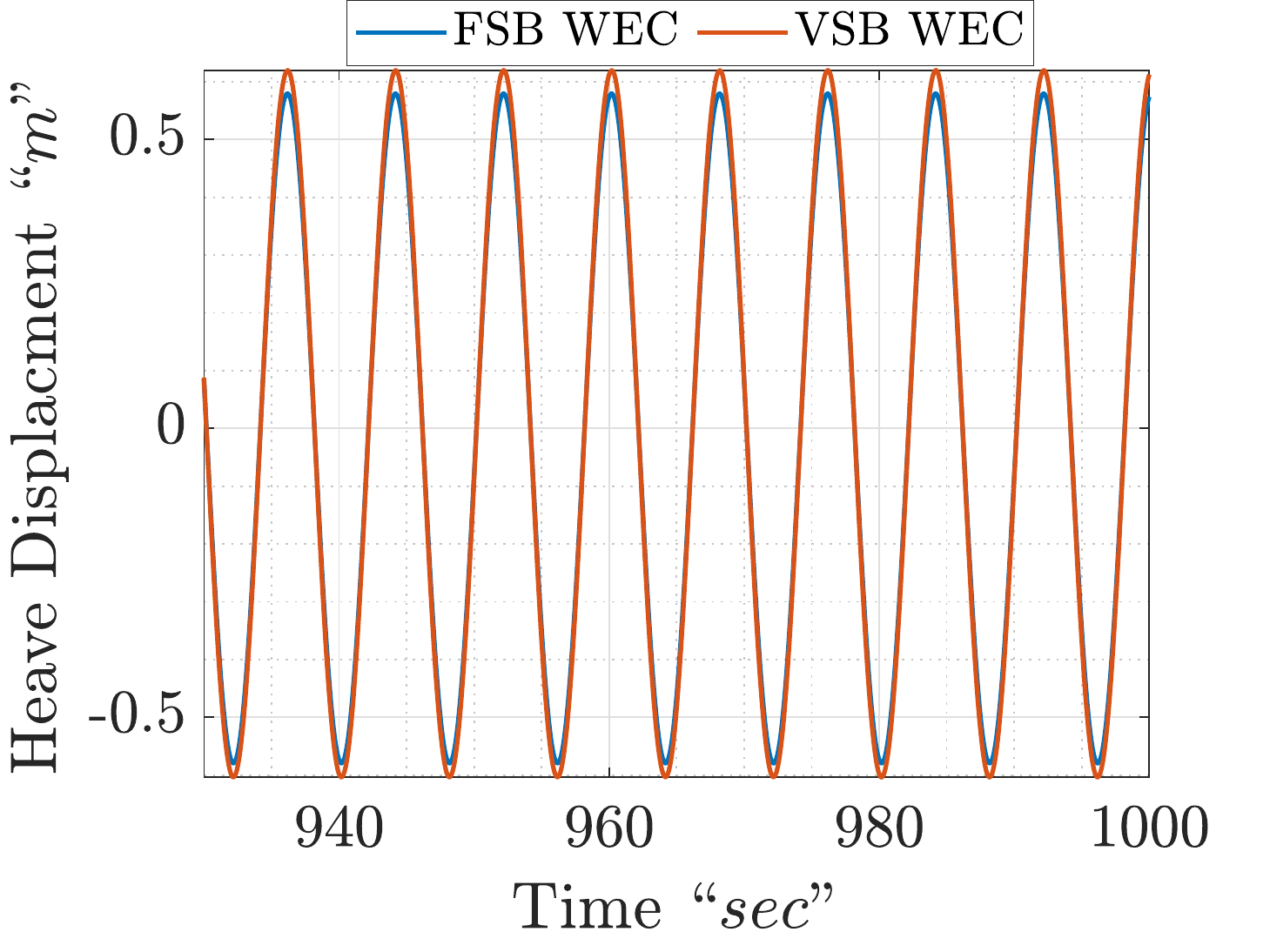}
\caption{Heave Displacement Response for Regular Wave (ID: RS14)}
\label{fig:heave_reg}
\end{figure}

Figure \ref{fig:heave_reg} shows a 30-second interval of the displacement response of the VSB WEC and the FSB WEC for the regular wave condition RS14. The results show an increase of 6.89 \% in the heave response for the VSB WEC compared to the FSB WEC.   

\begin{figure}%
    \centering
    \subfloat[\centering The Heave Velocity ]{{\includegraphics[scale=0.5]{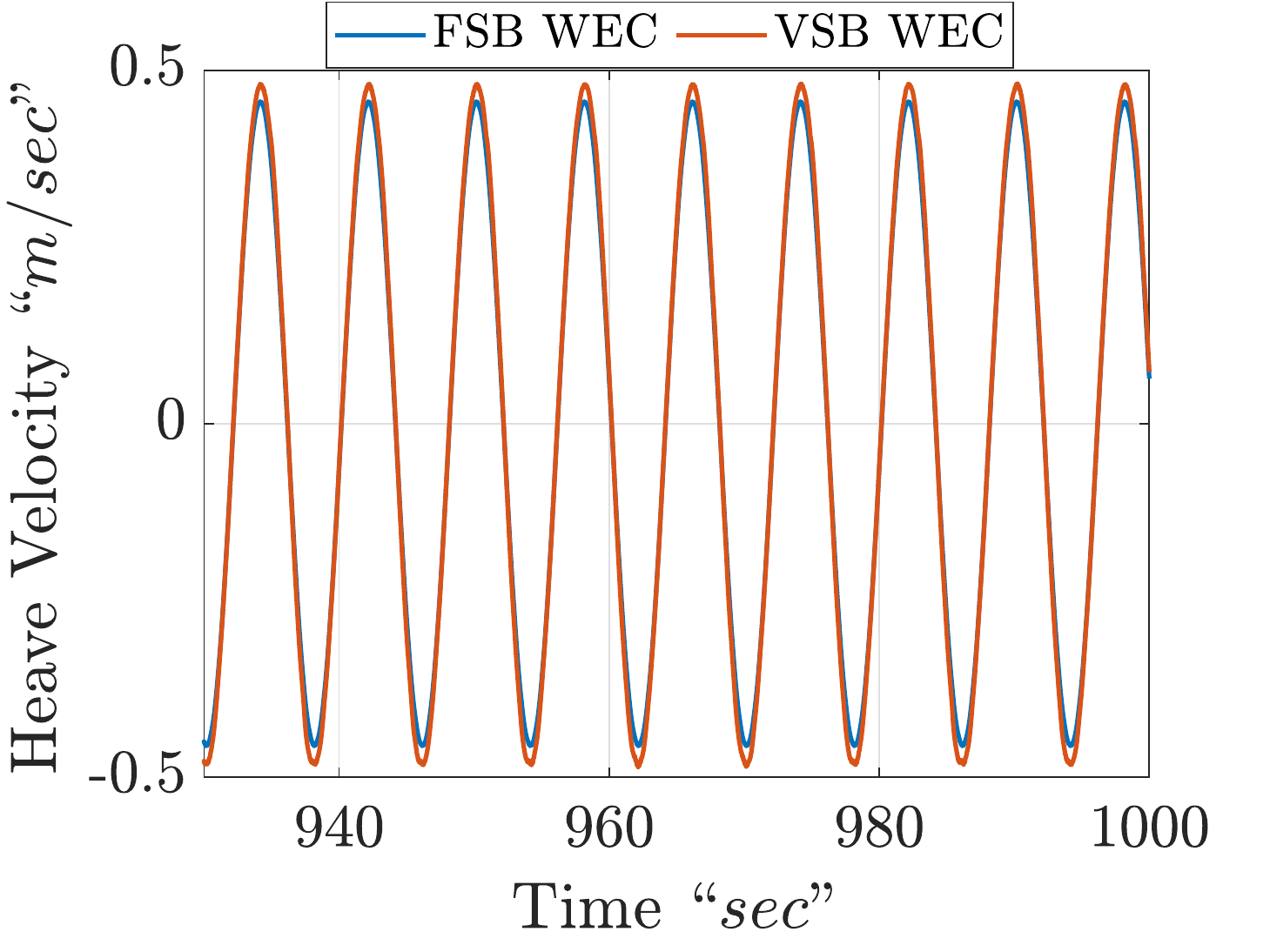}}} %
    \qquad
    \subfloat[\centering The PTO Force]{{\includegraphics[scale=0.5]{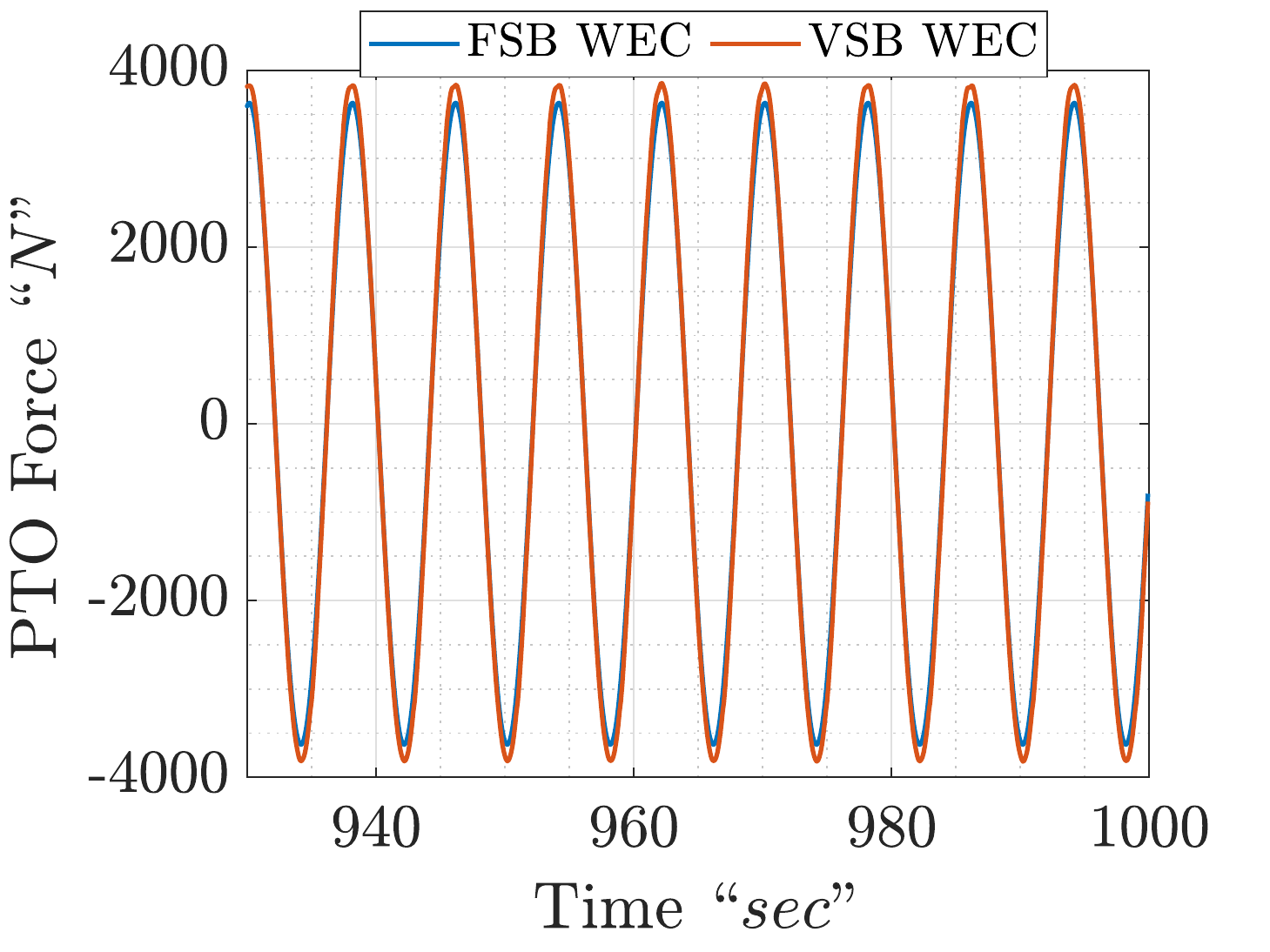} }} 
    \caption{Heaving Velocity and PTO Force for ID: RS14}%
    \label{fig:res2}%
\end{figure}

\begin{figure}%
    \centering
    \subfloat[\centering The Generated Power]{{\includegraphics[scale=0.5]{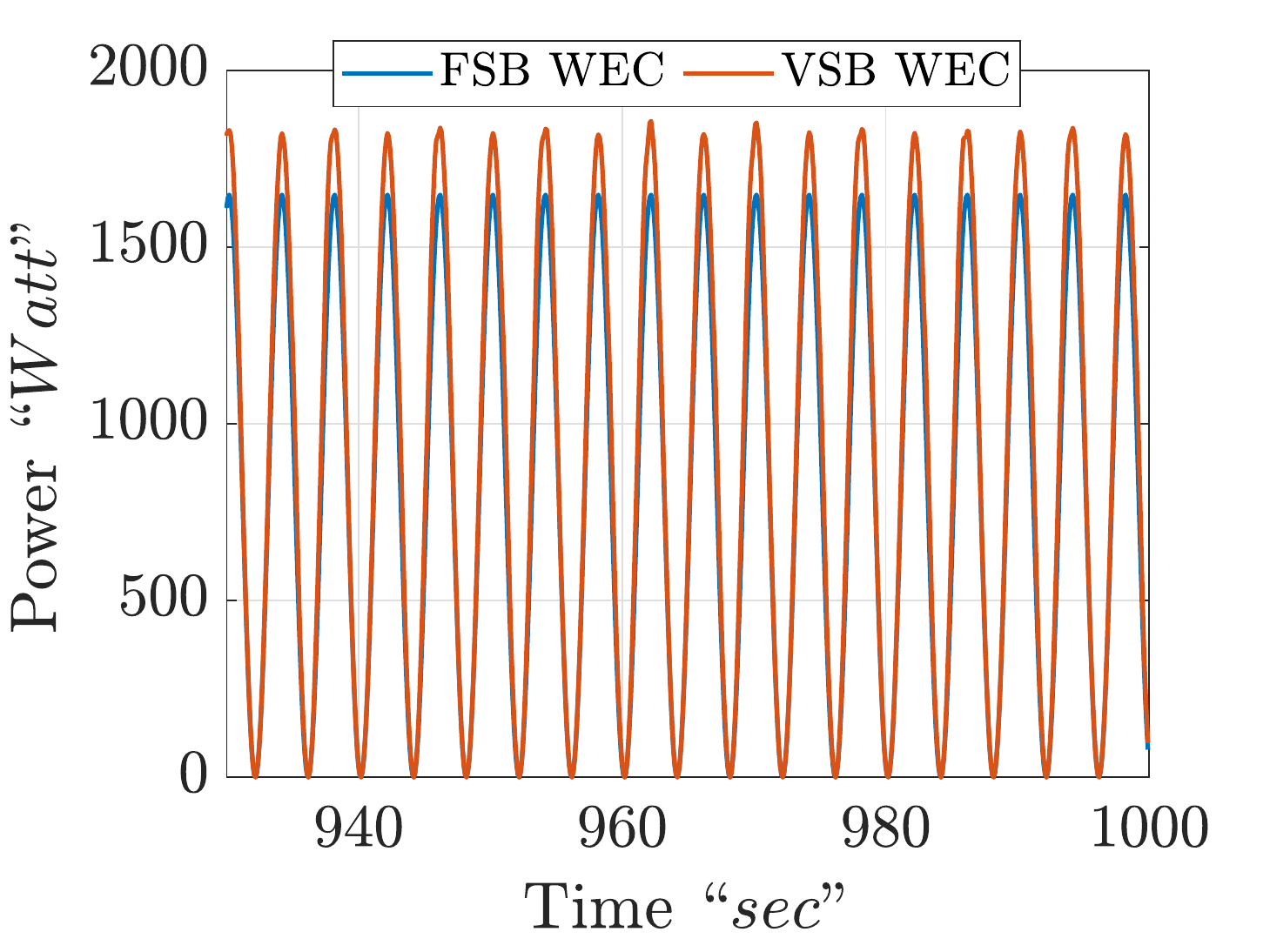}}} %
    \qquad
    \subfloat[\centering The Harvested Energy]{{\includegraphics[scale=0.5]{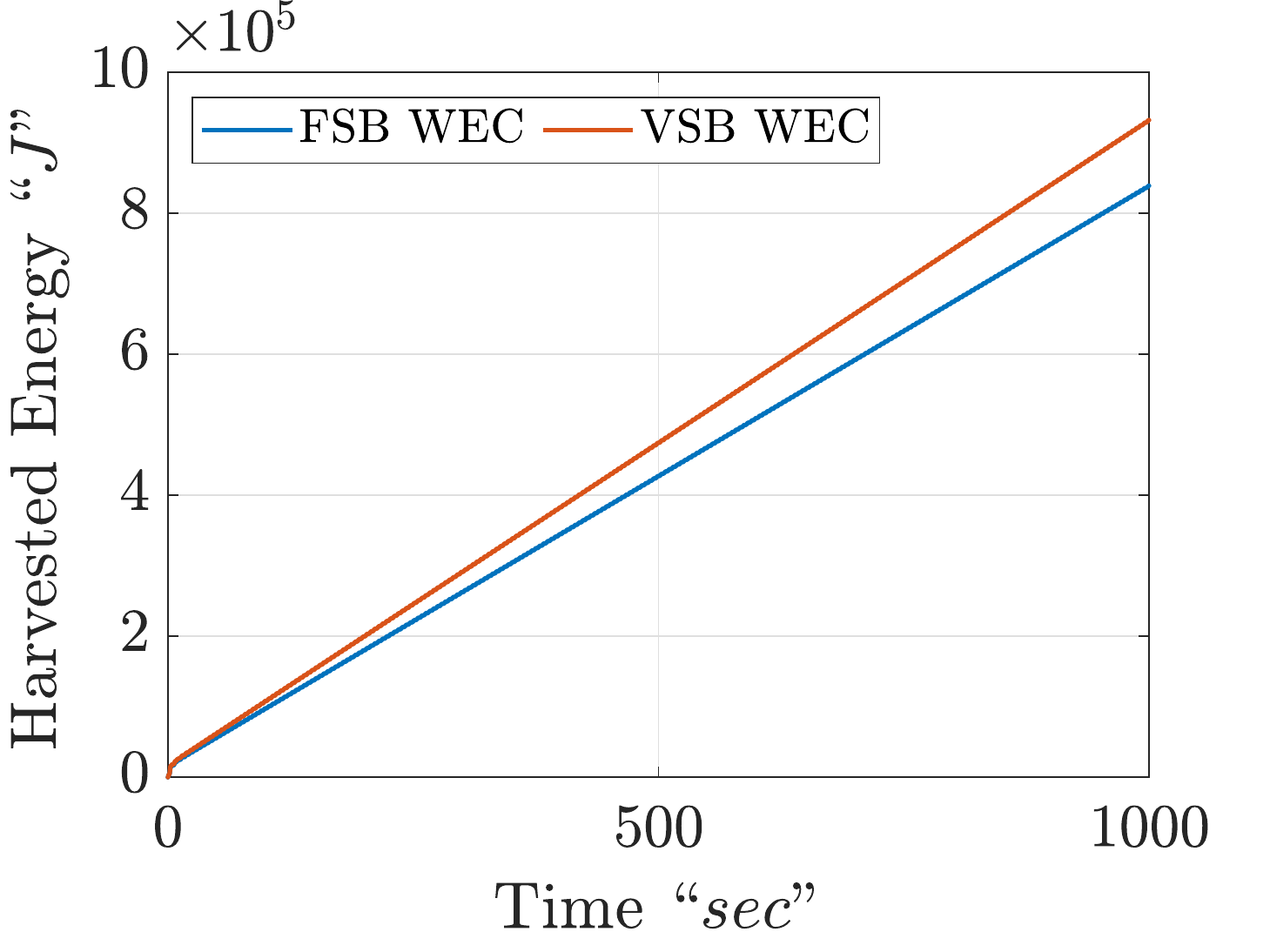} }} 
    \caption{Generated Power and Harvested Energy for ID: RS14}%
    \label{fig:respower}%
\end{figure}

The increase in heave response is accompanied by an increase of 5.58 \% in the velocity response for the VSB WEC, as shown in Fig.~\ref{fig:res2}-a. A similar increase in the PTO force is obtained but with $180^o$ phase shift between the velocity and force responses (Fig.~\ref{fig:res2}-b). Also, the power output increased by 12.7 \%, and the harvested energy increased by 11 \%  as shown in Fig.~\ref{fig:respower}.

Table~\ref{tab:3} shows a comprehensive comparison between the performance of the FSB WEC vs. the VSB WEC for the tested sea conditions; it can be seen that the VSB WEC harvested more energy in all the tested wave conditions with an increase in harvested energy ranging from 7.6 \% to 15.84 \%. It can be observed that the heave displacement and velocity pk-pk response increased for the VSB WECs over the FSB WEC for similar incident waves and PTO unit. This derives the conclusion that the VSB experiences more excited than the FSB for identical operating conditions. Also, larger passive controlled PTO units can be used for the same VSB WEC dimension.

\begin{table}[]
\caption{Regular Waves Results for the Wave Conditions and simulation time in Table \ref{tab:2}}\label{tab:3}

\begin{tabular}{|c|c|cc|cc|cc|c|}
\hline
\multirow{2}{*}{\textbf{ID}} & \multicolumn{1}{l|}{\multirow{2}{*}{$h$  (m)}} & \multicolumn{2}{c|}{Displacement (pk-pk m)} & \multicolumn{2}{c|}{Velocity (pk-pk m/s)} & \multicolumn{2}{c|}{Energy (J)}        & \multirow{2}{*}{Energy Increase} \\ \cline{3-8}
                             & \multicolumn{1}{l|}{}                          & \multicolumn{1}{c|}{FSB WEC}    & VSB WEC   & \multicolumn{1}{c|}{FSB WEC}   & VSB WEC  & \multicolumn{1}{c|}{FSB WEC} & VSB WEC &                                  \\ \hline
RS06                         & 0.1                                            & \multicolumn{1}{c|}{0.12548}    & 0.13498   & \multicolumn{1}{c|}{0.3168}    & 0.3398   & \multicolumn{1}{c|}{29920}   & 33960   & 13.5 \%                          \\ \hline
RS07                         & 0.1                                            & \multicolumn{1}{c|}{0.3054}     & 0.3178    & \multicolumn{1}{c|}{0.6422}    & 0.6636   & \multicolumn{1}{c|}{122300}  & 131600  & 7.6 \%                           \\ \hline
RS08                         & 0.1                                            & \multicolumn{1}{c|}{0.3986}     & 0.4168    & \multicolumn{1}{c|}{0.7142}    & 0.7438   & \multicolumn{1}{c|}{153000}  & 166600  & 8.9 \%                           \\ \hline
RS09                         & 0.15                                           & \multicolumn{1}{c|}{0.4814}     & 0.5194    & \multicolumn{1}{c|}{0.755}     & 0.81     & \multicolumn{1}{c|}{172400}  & 199700  & 15.84 \%                         \\ \hline
RS10                         & 0.2                                            & \multicolumn{1}{c|}{0.5672}     & 0.6002    & \multicolumn{1}{c|}{0.7892}    & 0.8318   & \multicolumn{1}{c|}{189200}  & 210900  & 11.47 \%                         \\ \hline
RS11                         & 0.2                                            & \multicolumn{1}{c|}{0.6532}     & 0.693     & \multicolumn{1}{c|}{0.819}     & 0.8636   & \multicolumn{1}{c|}{204600}  & 229200  & 12.02 \%                         \\ \hline
RS12                         & 0.2                                            & \multicolumn{1}{c|}{0.8246}     & 0.8874    & \multicolumn{1}{c|}{0.8614}    & 0.8952   & \multicolumn{1}{c|}{229500}  & 261100  & 13.77 \%                         \\ \hline
RS13                         & 0.25                                           & \multicolumn{1}{c|}{0.9938}     & 1.06      & \multicolumn{1}{c|}{0.889}     & 0.9254   & \multicolumn{1}{c|}{295500}  & 327000  & 10.66 \%                         \\ \hline
RS14                         & 0.25    & \multicolumn{1}{c|}{1.1588}           &      1.2386     & \multicolumn{1}{c|}{0.9108}          &  0.9616        & \multicolumn{1}{c|}{312100}     &    3412100     &    10.99    \%                          \\ \hline
\end{tabular}
\end{table}

\subsection{One-way vs. Two-way FSI results Comparison:}
The water free-surface location does not change during the simulation, which is also the mean value of the incident waves as the time approaches infinity, i.e., the buoy C.G location oscillates about the free surface. The VSB in the simulations has a perfectly spherical shape at $t=0$ (similar to the FSB). When the waves excite the vibration modes of the sphere, it can be noticed that each one of the mode shapes vibrates around the initial shape of the VSB. 

It can also be noticed when plotting the variation of each element in the generalized added mass, generalized hydrodynamic damping, and the generalized hydrostatic matrices that they all oscillate around their initial values at the beginning of the simulations. Similar observations can be made when plotting the flexibility states ($\boldsymbol{\eta}$). Accordingly, the mean values of these matrices are equal to their corresponding matrices for the FSB.\\

% \usepackage{multirow}
\begin{table}[]
\caption{Comparison between the simulation time and errors for one-way and two-way FSI for the VSB WEC models}\label{tab:4}
\begin{tabular}{|c|cc|c|c|}
\hline
\multirow{2}{*}{\textbf{ID}} & \multicolumn{1}{c|}{One-Way FSI} & Two-Way FSI & \multicolumn{1}{l|}{FSI Simulation Time} & Energy Harvesting Error   \\ \cline{2-3}
                             & \multicolumn{2}{c|}{Computational Time (hour)} & (sec)                                    & \multicolumn{1}{c|}{(\%)} \\ \hline 
RS06        & \multicolumn{1}{c|}{17 (sec)}     &    39.42  (hour)   & 300   &  -1.98  \\ \hline
RS07     & \multicolumn{1}{c|}{0.19}        &    39.72  (hour) & 300  &    0.3    \\ \hline
RS08     & \multicolumn{1}{c|}{0.18}        &    40.69  (hour) & 300  &   -0.17  \\ \hline
RS09     & \multicolumn{1}{c|}{0.28}           &    38.64   (hour)   & 300  &  1.49      \\ \hline
RS10     & \multicolumn{1}{c|}{0.34}           &    37.85  (hour) & 300    & -0.23     \\ \hline
RS11        & \multicolumn{1}{c|}{0.20}        &    36.64  (hour) & 300    &  -0.15 \\ \hline
RS12     & \multicolumn{1}{c|}{0.21}           &    37.71  (hour)  & 300   &    0.94   \\ \hline
RS13        & \multicolumn{1}{c|}{0.2}        &    41.6   (hour)    & 360   &  0.77  \\ \hline
RS14        & \multicolumn{1}{c|}{0.21}        &    50   (hour)  & 360     & 0.95 \\ \hline
\end{tabular}
\end{table}

Table.~\ref{tab:4} shows the results obtained for regular waves from one-way and two-way FSI models. The one-way FSI model used in this work uses averaged hydrodynamic coefficients for the VSB WEC discussed in subsection \ref{sec:one-two_way}. Running one-way FSI reduced the computation time from days to a few minutes for regular waves, and from weeks to a couple of days for irregular waves. Table~\ref{tab:4} shows a comparison between the energies and the errors between the two models for the case of regular waves. It can be seen that the reduction in simulation time is significant, and the error in harvested energy is bounded between $\pm 1.98$. 

The simulation of two-way FSI for irregular waves is computationally expensive as it requires weeks to get one simulation done. On the other hand, the one-way FSI simulation for irregular waves took a couple of days to get an approximately similar result. For this reason, all the irregular waves simulations in the next subsection were done using the one-way FSI model. 

\subsection{Irregular Waves}
The Bretchsneider wave spectrum was used to generate irregular waves with $N_\omega =256$ equidistant frequencies ranging between 0.1 to 3.5 rad/sec. The one-way FSI Eq.~\eqref{eq:EOM_1w_2} was applied. The shell thickness was tuned so that the shell deformation remained within a certain bound: A large thickness makes the shell rigid, and a small thickness would yield unrealistic deformations. This is because the shell strain expressions in Eq.~\eqref{eq:15.1} and \eqref{eq:16.1} are based on a small deformations assumption \cite{shabara2021numerical}.

\begin{figure}[h]
\centering
\includegraphics[scale=0.6]{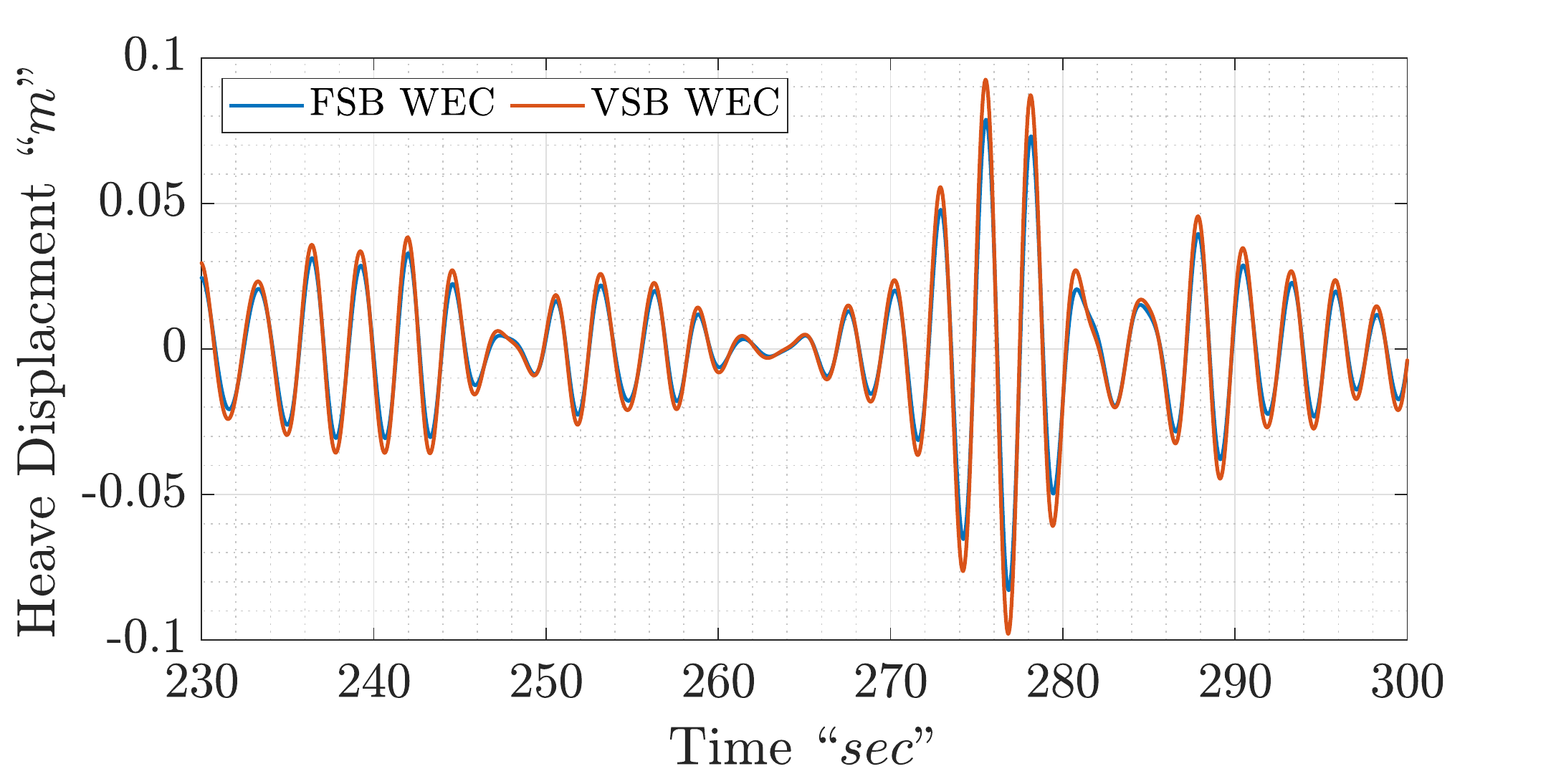}
\caption{Steady State Heave Displacement Response for Irregular Wave (ID: RS06)}
\label{fig:irr_1}
\end{figure}

\begin{figure}%
    \centering
    \subfloat[\centering Heave velocity response]{{\includegraphics[scale=0.5]{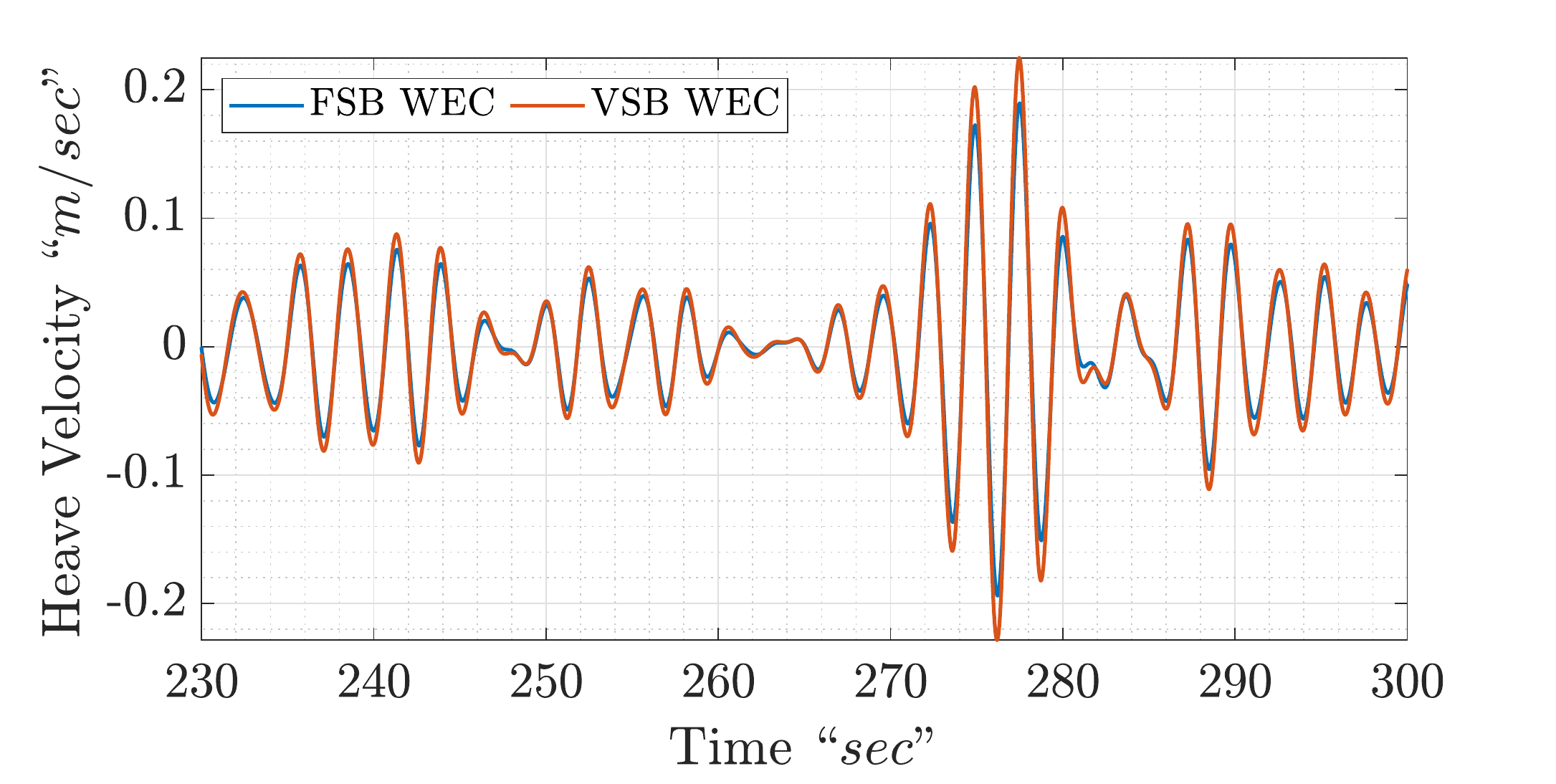}}} %
    \qquad
    \subfloat[\centering Power take-off Unit Force Energy]{{\includegraphics[scale=0.5]{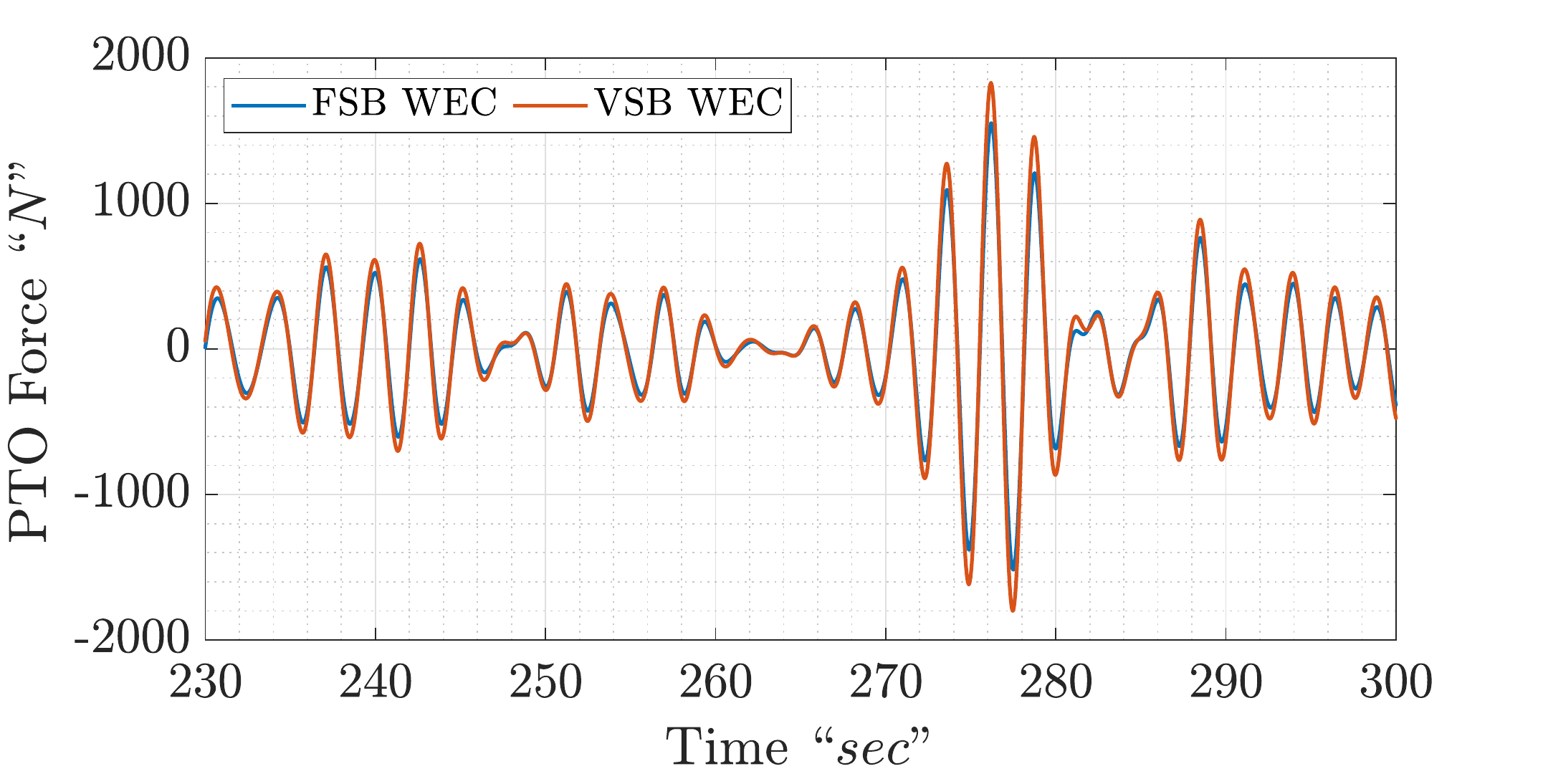} }} 
    \caption{Heave velocity response and PTO force waveform for irregular waves (ID: RS06)}%
    \label{fig:irr_2}%
\end{figure}

The harmonic response of the VSB WEC decays within the first 10 seconds of the simulation. Fig.~\ref{fig:irr_1} shows an interval of 70 seconds of the heave response of both FSB and VSB WECs in response to the RS06 irregular wave. The pk-pk heave displacement increased by a factor of 18\% for the VSB WEC. Similarly, the pk-pk heave velocity and PTO responses increased by a factor of 18\%, as shown in Fig.~\ref{fig:irr_2}.

\begin{figure}[h]
\centering
\includegraphics[scale=0.6]{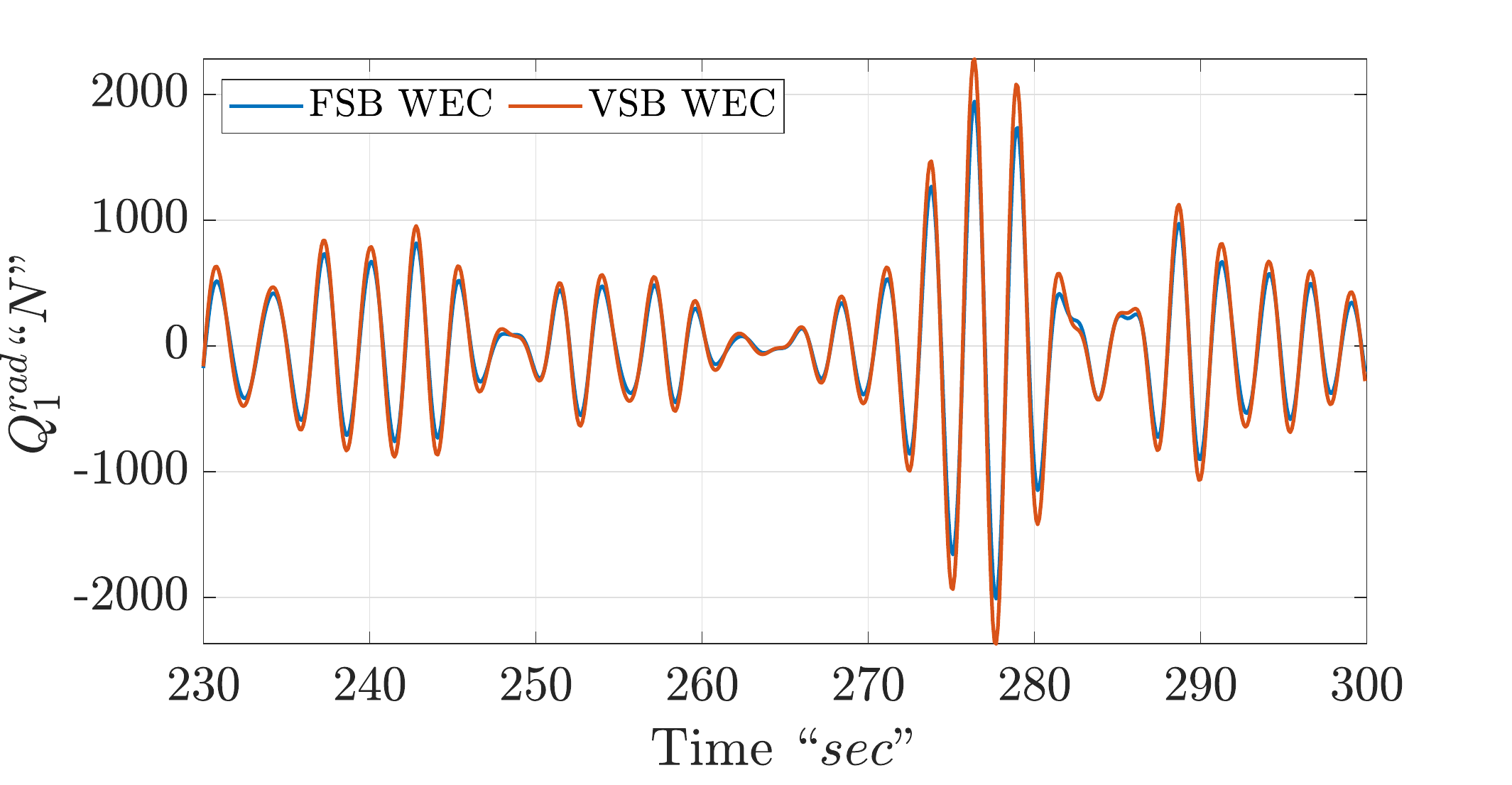}
\caption{Generalized radiation force $\boldsymbol{Q}^{rad}_1$ on buoys' CG (ID: RS06)}
\label{fig:irr_3}
\end{figure}

\begin{figure}[h]
    \centering
    \subfloat[\centering The Generated Power]{{\includegraphics[scale=0.5]{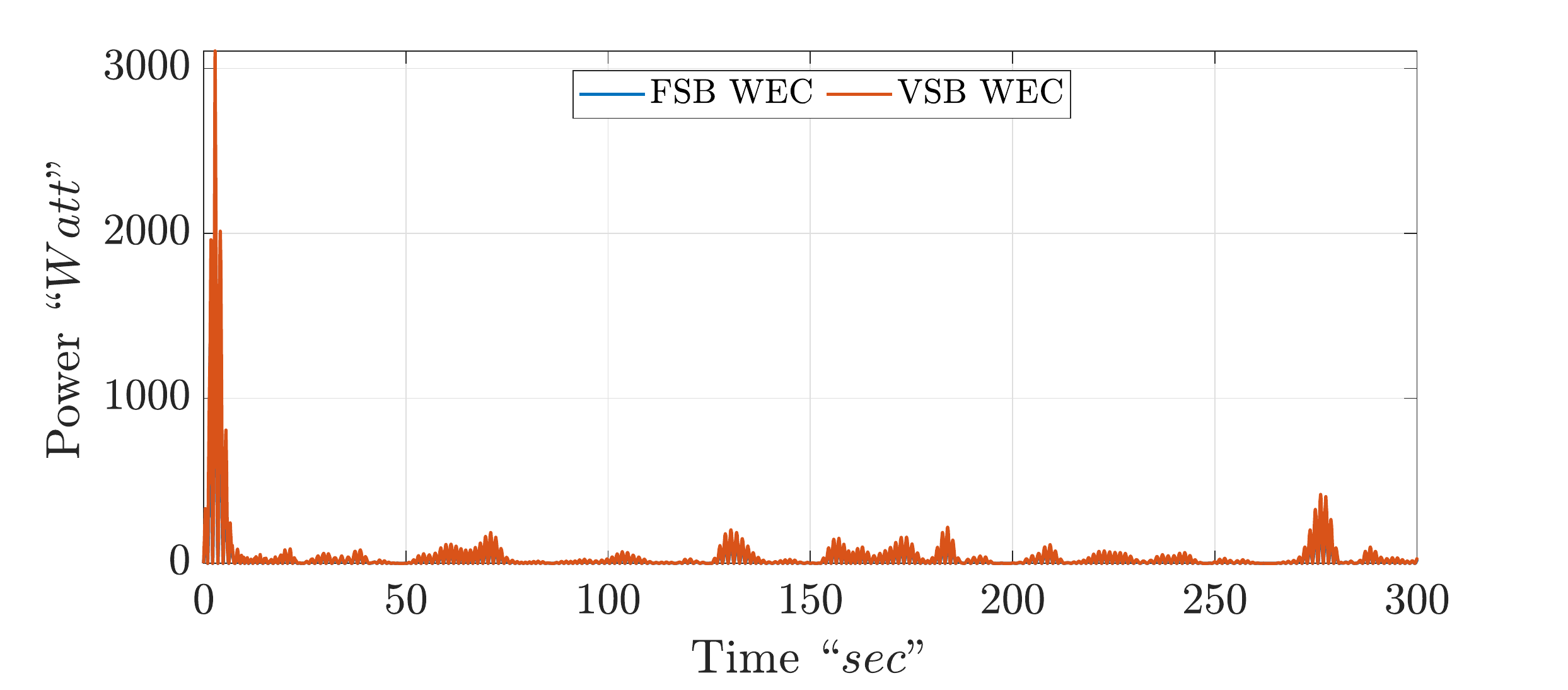}}} %
    \qquad
    \subfloat[\centering The Harvested Energy]{{\includegraphics[scale=0.5]{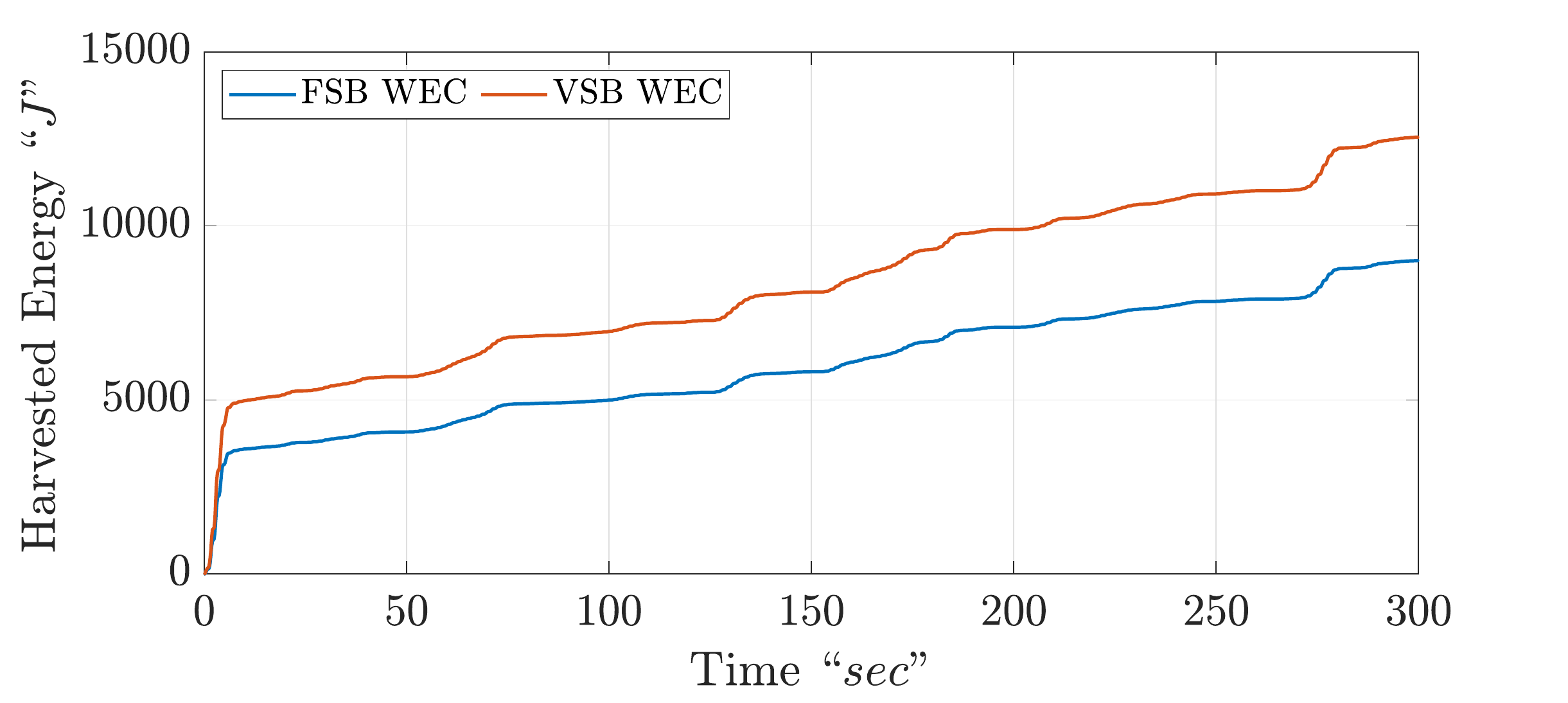} }} 
    \caption{Generated Power and Harvested Energy for Irregular Waves ID: RS06}%
    \label{fig:irr_4}%
\end{figure}

Figure \ref{fig:irr_3} shows a plot of the first component of the generalized radiation force (acting on the buoy's C.G)  $\boldsymbol{Q}^{rad}_1$, which is calculated using Eq.~\eqref{eq:55_2}. It can be noticed that the pk-pk radiation force increased when using VSB. Fig.~\ref{fig:irr_4}-a shows the full simulation interval for the generated powers; it can be noticed that the largest power peaks are generated in the transient response region (first 10 seconds). This is also demonstrated in harvested energy presented in Fig.~\ref{fig:irr_4}-b. This agrees with the results obtained from high-fidelity simulations (FEA and CFD) in \cite{shabara2021numerical}.

Table \ref{tab:5} shows the results between the FSB and VSB  WECs for the tested wave conditions. It can be noticed that there is an increase of almost 40\% in the harvested energy harvested by the VSB WEC for the RS06 wave condition. The other wave conditions showed an increase in harvested energy with different degrees. The main RS06 model used $E = 1e5$Pa compared to $2e5$Pa for the other tested wave conditions. Also, smaller thicknesses result in more energy harvesting as it is translated to more flexibility for the shell, which results in more excitation of the shell.

\begin{table}[]
\caption{Irregular Waves Results for the Wave Conditions and simulation time in Table \ref{tab:2}}\label{tab:5}

\begin{tabular}{|c|c|cc|cc|cc|c|}
\hline
\multirow{2}{*}{\textbf{ID}} & \multicolumn{1}{l|}{\multirow{2}{*}{$h$  (m)}} & \multicolumn{2}{c|}{Displacement (pk-pk m)} & \multicolumn{2}{c|}{Velocity (pk-pk m/s)} & \multicolumn{2}{c|}{Energy (J)}        & \multirow{2}{*}{Energy Increase} \\ \cline{3-8}
                             & \multicolumn{1}{l|}{}                          & \multicolumn{1}{c|}{FSB WEC}    & VSB WEC   & \multicolumn{1}{c|}{FSB WEC}   & VSB WEC  & \multicolumn{1}{c|}{FSB WEC} & VSB WEC &                                  \\ \hline
RS06  & 0.11    & \multicolumn{1}{c|}{0.16}  &  0.19&\multicolumn{1}{c|}{0.3825}    &  0.4521  & \multicolumn{1}{c|}{9005}   &  12550  & 39.37\%                          \\ \hline
RS07  & 0.1    & \multicolumn{1}{c|}{0.287}  &  0.3012 & \multicolumn{1}{c|}{0.6443}    &  0.6761  & \multicolumn{1}{c|}{27340}   &  30090  & 10\%                          \\ \hline
RS08   & 0.1       & \multicolumn{1}{c|}{0.409}     &  0.4302   &  \multicolumn{1}{c|}{0.8623}    & 0.9064   & \multicolumn{1}{c|}{55580}  & 61390  & 10.45 \%                   \\ \hline
RS09   &   0.15    & \multicolumn{1}{c|}{0.5059}     & 0.5188    & \multicolumn{1}{c|}{1.009}     &   1.0332   & \multicolumn{1}{c|}{91330}  & 95800  &  4.8\%                         \\ \hline
RS10   &   0.2    & \multicolumn{1}{c|}{0.5557}     & 0.562    & \multicolumn{1}{c|}{1.0727}     &   1.0737  & \multicolumn{1}{c|}{131200}  & 133362  &  1.8\%                         \\ \hline
RS11     &   0.2      & \multicolumn{1}{c|}{0.6259}     & 0.6355    & \multicolumn{1}{c|}{1.05}     & 1.064   & \multicolumn{1}{c|}{168100}  &  172500 &  2.62\%                         \\ \hline
RS13     &  0.22   & \multicolumn{1}{c|}{0.9879}           &   1.002        & \multicolumn{1}{c|}{1.196}          &  1.2131   & \multicolumn{1}{c|}{284200}     &    290200     &   2.2\%                          \\ \hline
\end{tabular}
\end{table}

The simulation time for irregular waves in a one-way FSI environment ranged from 30-40 hours which is a significant simulation time improvement over both two-way FSI models (low fidelity and the high fidelity \cite{shabara2021numerical}). 

\subsection{Spherical VSB Physics:}

\begin{figure}[h]
    \centering
    \subfloat[\centering Breathing Mode $n=0$ associated with $\eta_1$]{{\includegraphics[scale=0.45]{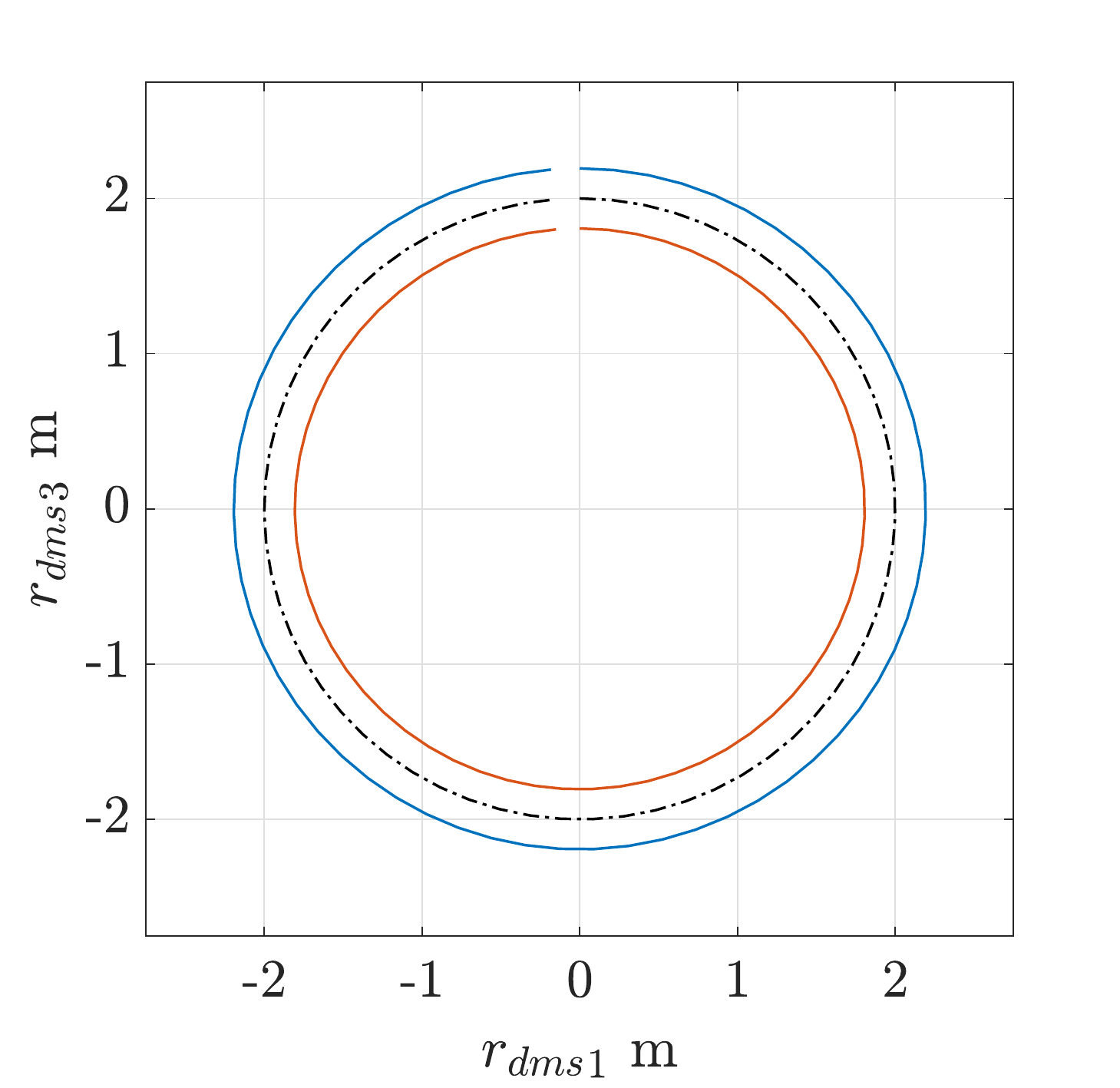}}} %
    \qquad
    \subfloat[\centering $\eta_2$]{{\includegraphics[scale=0.59]{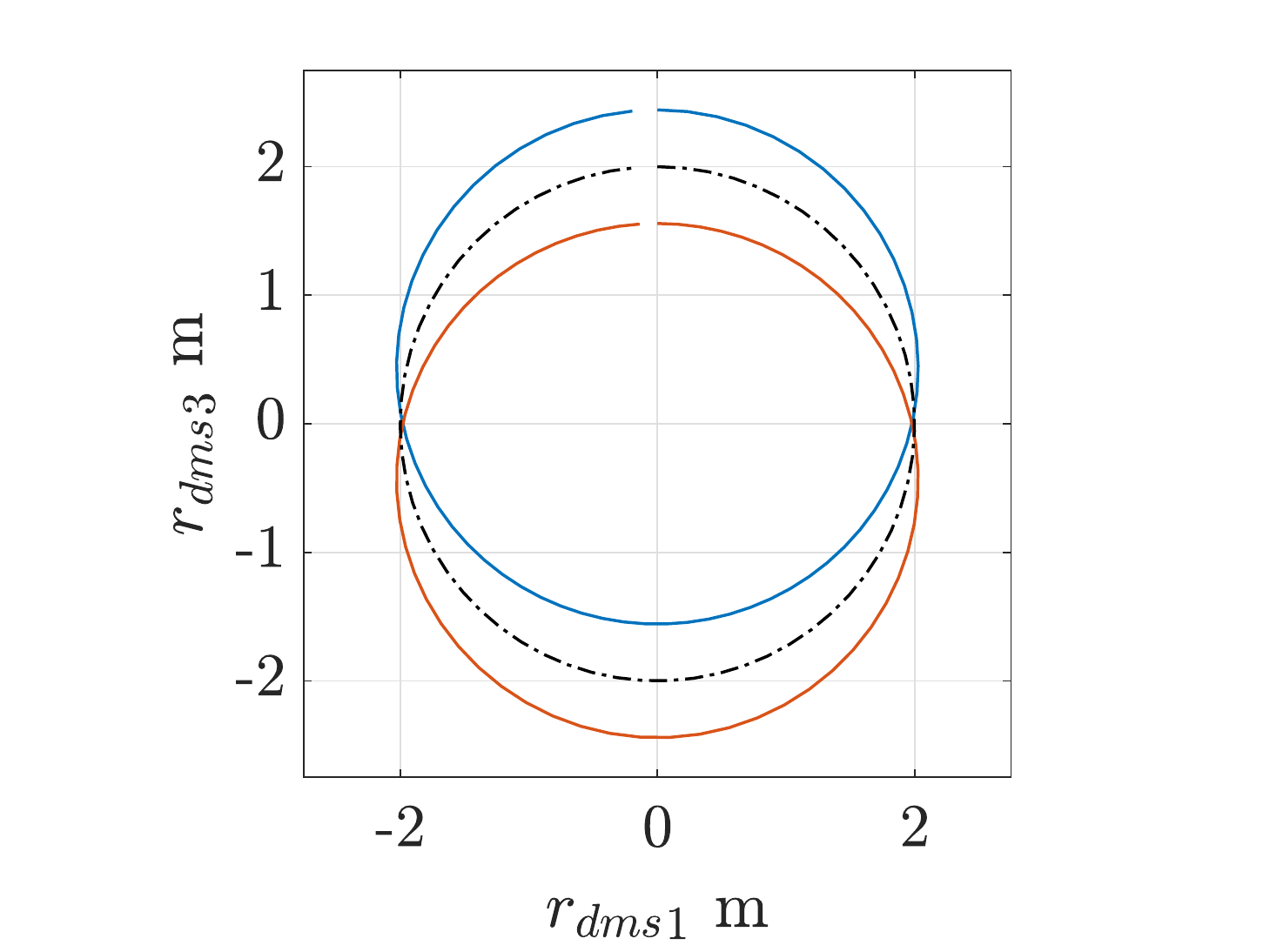} }} 
        \qquad
    \subfloat[\centering $\eta_3$]{{\includegraphics[scale=0.45]{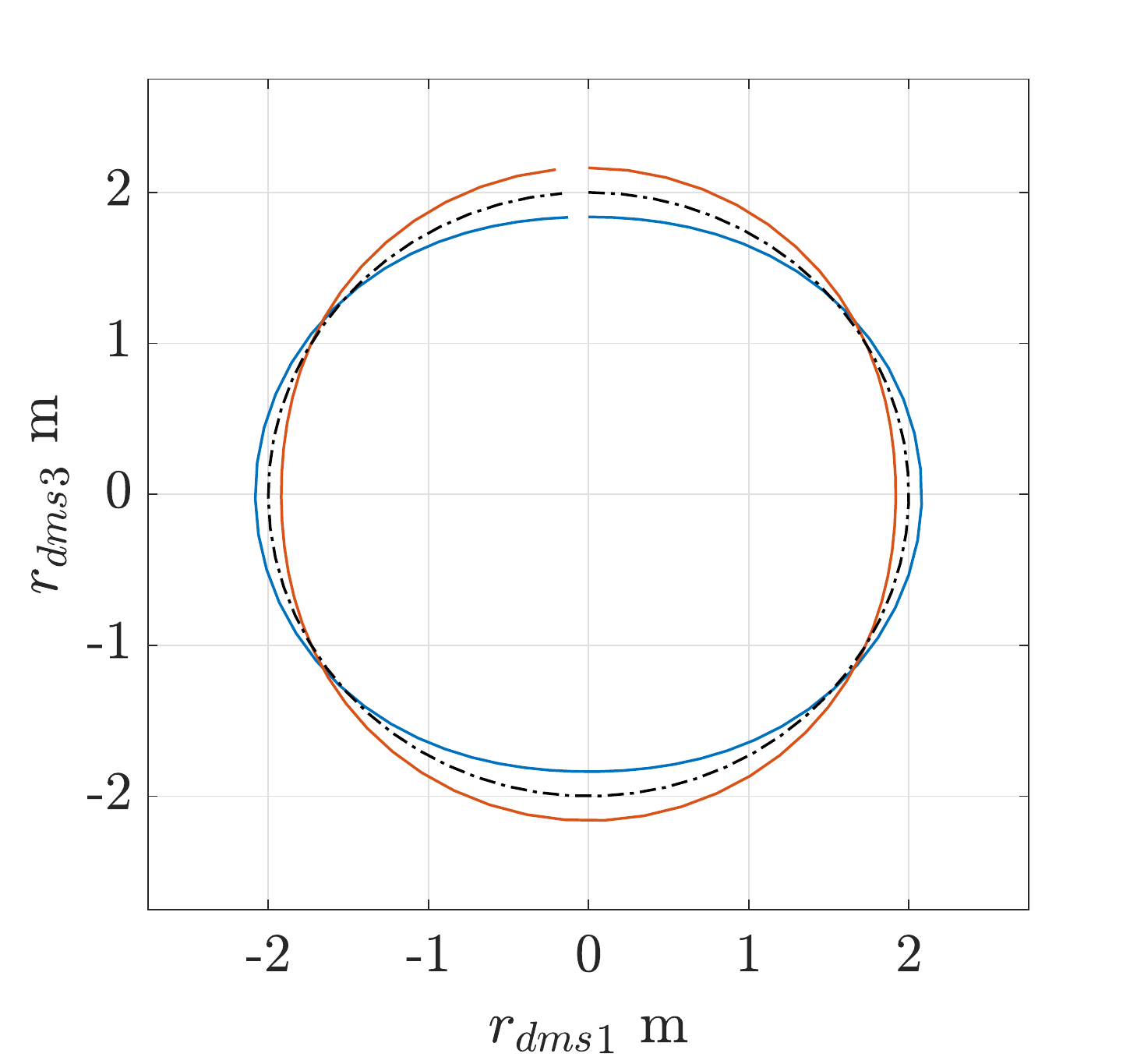} }} 
        \qquad
    \subfloat[\centering $\eta_4$]{{\includegraphics[scale=0.59]{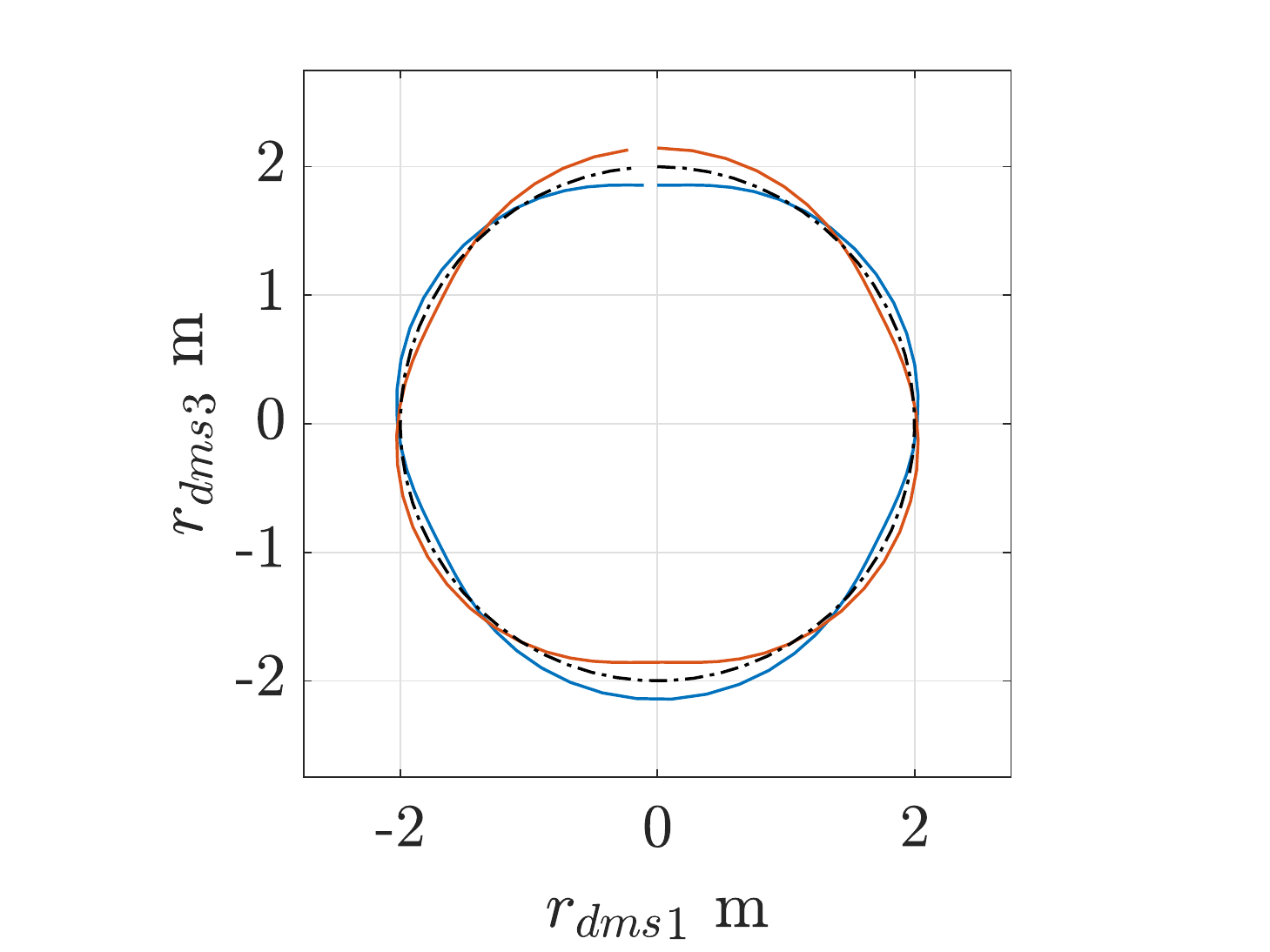} }} 
    \caption{Axisymmetrical vibration modes of a spherical shell, the black dotted line is the non-deformed shape and the colored lines are for the deformed shape}%
    \label{fig:modes}%
\end{figure}

The first four vibration modes of the spherical buoy are shown in Fig.~\ref{fig:modes}. The first mode shape corresponds to the flexibly state $\eta_1$ and the breathing mode $n=0$ where the shell is deformed uniformly in the radial direction. The second mode shape corresponds to the flexibility state $\eta_2$ and $n=1$, where the shell encounters the rigid-body motion, and the buoy C.G encounters pure translational motion without any shell deformations. The third mode shape corresponding to $n=1$ and $\eta_3$ deform the shell to become either prolate spheroid when $\eta_3$ is positive, or oblate spheroid when $\eta_3$ is negative. Noting that, the second and third mode shapes both corresponds to $n=1$ but with the change in sign in Eq.~\eqref{eq:7}. Finally, the fourth mode shape corresponds to $\eta_4$, resulting in a deformed shape that switches from pointy to plat shapes. The dominant mode shapes throughout the simulations were two of the first three mode shapes, and their activation contributed to the increase in vibration responses and harvested energy for the VSB WECs over the FSB WECs.

\begin{figure}[h]
\centering
\includegraphics[scale=0.45]{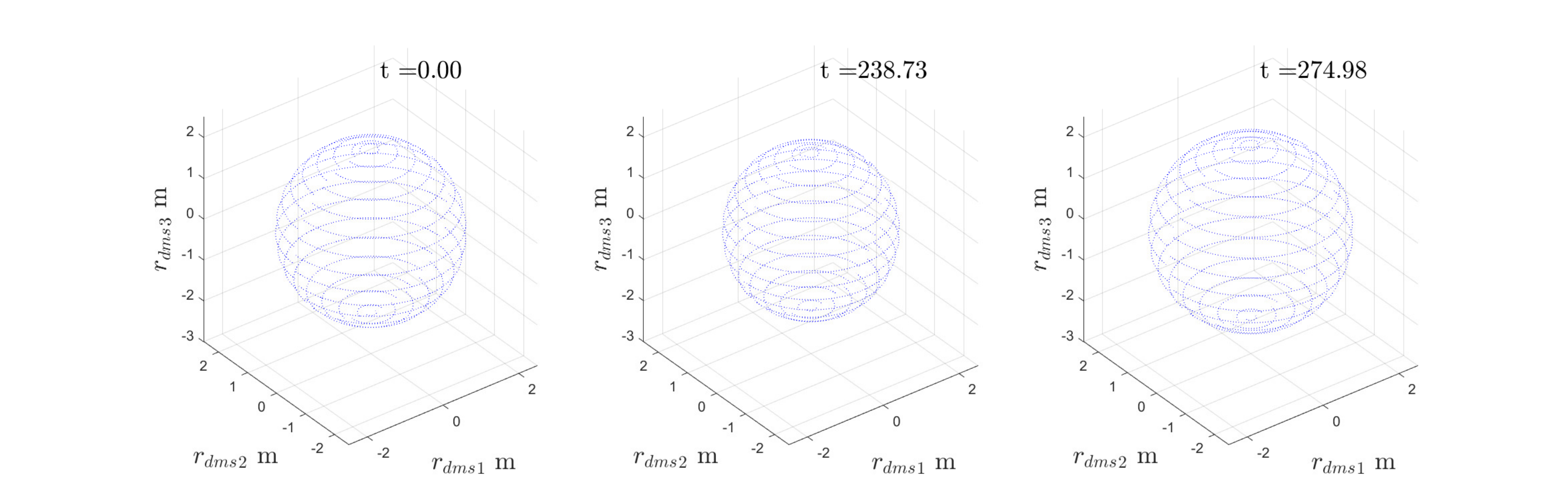}
\caption{3D plot for VSB WEC at different time instances (Deformation vector scale 5:1)  (ID: RS06)}
\label{fig:3dwec}
\end{figure}

Figure \ref{fig:3dwec} shows the 3D plot for the VSB WEC at different time instances for wave condition RS06. At $t=0$ the shell has a pure spherical shape (initially non-deformed). At $t=238.73$ the buoy is encountering a beginning of a wave trough, and at $t=274.98$ is encountering a wave crest. 
The first and the third mode shapes are dominant for the RS06 wave condition. The beginning of the wave crest coincides with the shell expansion related to the first mode and the oblate spheroid shape related to the third mode shape; this increases the vertical component of the excitation force on the C.G and the heave excitation in the positive vertical direction. On the other hand, the contraction deformation and the prolate spheroid shapes decrease the vertical component of the excitation force of the buoy C.G; this causes the buoy to plunge deeper into the waves when it encounters a trough compared to the FSB WEC. The excitation of these mode shapes causes an increase in the harvested energy.

% \usepackage{multirow}
% \usepackage{multirow}

\clearpage

\section{Conclusion and Future Work} \label{sec:conclusion}
This paper derives a novel equation of motion for variable shape wave energy converters; this equation of motion can be viewed as a modified Cummins equation for flexible buoys. The proposed equation is derived in the context of Lagrangian mechanics. Novel expressions were derived for the generalized hydrodynamic forces/coefficients (added mass, damping, hydrostatic and excitation) that account for the buoy flexibility in regular and irregular waves. To reduce the computational cost accompanied with the proposed two-way Fluid-Structure Interaction scheme, the Reynolds averaging technique was used to derive average hydrodynamic coefficients. These average hydrodynamic coefficients enable one-way FSI model. The model was tested using spherical flexible shape WECs. The numerical results support this work's hypothesis that a VSB WEC would harvest energy at a higher rate than that of an FSB WEC when both WECs use no reactive power. 

For future work, optimal control methods have to be developed to control the power take-off unit, and the shape of the buoy as a function of the incident wave conditions. Also, the analytical model can be extended to model hyperelastic materials, smart materials, power take-off unit dynamics, and other WEC geometries. 

\section*{Acknowledgments}
This material is based upon work supported by the National Science Foundation (NSF), USA, under Grant Number 2023436.
The research reported in this paper is partially supported by the HPC@ISU equipment at Iowa State University, some of which has been purchased through funding provided by NSF under MRI grant number 1726447.


%\bibliographystyle{ieeetr}
%\printbibliography
%\bibliographystyle{natbib}
\bibliographystyle{unsrtnat}

%\bibliographystyle{cas-model2-names}
%\bibliographystyle{abbrvnat}

\bibliography{sample,myRefs,ORefs,WEC}

\bio{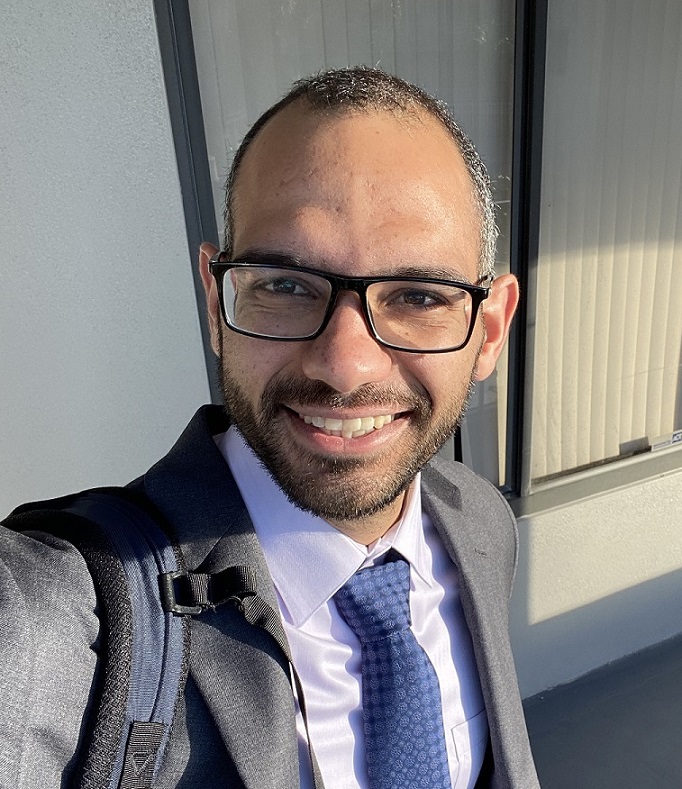}
Mohamed Shabara is a Ph.D. Student in the aerospace department at Iowa State University. He received his M.S. degree in Mechanical and Aerospace engineering from Old Dominion University. Also, he received M.S. and B.S. in mechanical engineering from Alexandria University in Egypt. Mr. Mohamed worked as a turbo-machinery engineer in the Oil and Gas fields. His research interest includes dynamics and optimal control
\endbio  

\vspace{1cm}

\bio{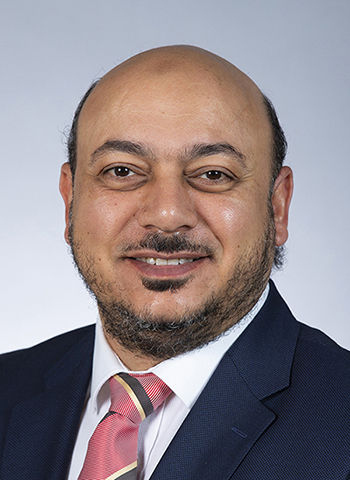}
Ossama Abdelkhalik is a Professor at Iowa State University. He received his PhD degree  in Aerospace Engineering from Texas A\&M University, and his BS and MS degrees from Cairo University in Egypt. Dr. Abdelkhalik’s research is in the area of global optimization, dynamics, and control with  applications to the wave energy conversion and space trajectory optimization. Dr. Abdelkhalik is Associate Fellow of AIAA and is associate editor of the AIAA Journal of Spacecraft and Rockets.
\endbio

\end{document}